\documentclass{article}
\usepackage{amsmath,amsfonts,amsthm,amssymb,amscd,bm,bbm}
\usepackage[hidelinks]{hyperref}
\usepackage[usenames,dvipsnames]{color}
\usepackage{url}
\hypersetup{colorlinks, citecolor=ForestGreen, linkcolor=MidnightBlue, urlcolor=Blue}

\newcommand{\be}{\begin{equation}}
\newcommand{\ee}{\end{equation}}
\newcommand{\bp}{\begin{proof}}
\newcommand{\ep}{\end{proof}}
\newcommand{\bi}{\begin{itemize}}
\newcommand{\ei}{\end{itemize}}
\newcommand{\om}{\omega}

\newcommand{\iI}{\mathfrak{i}}
\newcommand{\pP}{\mathfrak{p}}
\newcommand{\uu}{\mathfrak{u}}
\newcommand{\vv}{\mathfrak{v}}
\newcommand{\Zz}{\mathfrak{z}}
\newcommand{\aaa}{{\cal A}}
\newcommand{\BBB}{\mathfrak B}
\newcommand{\E}{\mathfrak E}

\newcommand{\ees}{{\cal E}}
\newcommand{\fff}{{\cal F}}

\newcommand{\gi}{{\cal G}}
\newcommand{\h}{{\cal H}}

\newcommand{\elll}{{\cal L}}
\newcommand{\mmm}{{\cal M}}

\newcommand{\ooo}{{\cal O}}

\newcommand{\vvv}{{\cal V}}

\newcommand{\yyy}{{\cal Y}}
\newcommand{\zzz}{{\cal Z}}

\newcommand{\e}{\mathbb{E}}

\newcommand{\nn}{\mathbb{N}}
\newcommand{\pp}{\mathbb{P}}

\newcommand{\rr}{\mathbb{R}}

\newcommand{\kp}{\varkappa}

\newcommand{\q}{\quad}

\newcommand{\f}{\frac}
\newcommand{\lm}{\lambda}
\newcommand{\p}{\partial}
\newcommand{\ph}{\varphi}

\newcommand{\De}{\delta}
\newcommand{\de}{\Delta}
\newcommand{\g}{\nabla}
\newcommand{\dt}{\dot}

\newcommand{\es}{\varepsilon}

\newcommand{\al}{\alpha}

\newcommand{\nt}{\noindent}

\newcommand{\ch} {\mathbbm{1}}

\newcommand{\iin} {\infty}

\newcommand{\mpp}{\emph}
\newcommand{\ef}{\eqref}
\newcommand{\dd}{\,{\textup d}}
\newcommand{\Dd}{{\textup d}}
\theoremstyle{plain}
\newtheorem{theorem}{Theorem}[section]
\newtheorem*{mt}{Main Theorem}
\newtheorem{lemma}[theorem]{Lemma}
\newtheorem{proposition}[theorem]{Proposition}

\theoremstyle{definition}
\newtheorem{definition}[theorem]{Definition}

\theoremstyle{remark}

\numberwithin{equation}{section}
\newtheorem*{definition*}{Definition}
\newtheorem*{problem*}{Problem}

\newtheorem*{remark*}{Remark}
\newtheorem*{note*}{Note}
\begin{document}
\author{Davit Martirosyan\footnote{Department of Mathematics, University of Cergy-Pontoise, CNRS UMR 8088, 2 avenue
Adolphe Chauvin, 95300 Cergy-Pontoise, France;e-mail: \href{mailto:Davit.Martirosyan@u-cergy.fr}{Davit.Martirosyan@u-cergy.fr}}}
\title{Large deviations for stationary measures of stochastic nonlinear wave equation with smooth white noise} 
\date{}
\maketitle

\begin{abstract}
The paper is devoted to the derivation of large deviations principle for the family $(\mu^\es)_{\es>0}$ of stationary measures of the Markov process generated by the flow of equation
$$
\p_t^2u+\gamma \p_tu-\de u+f(u)=h(x)+\sqrt{\es}\,\vartheta(t,x).
$$
The equation is considered in a bounded domain $D\subset\rr^3$ with a smooth boundary and is supplemented with the Dirichlet boundary condition. Here $f$ is a nonlinear term satisfying some standard dissipativity and growth conditions, the force $\vartheta$ is a non-degenerate white noise, and $h$ is a function in $H^1_0(D)$. The main novelty here is that we do not assume that the limiting equation (i.e., when $\es=0)$ possesses a unique equilibrium and that we do not impose roughness on the noise. Our proof is based on a development of the approach introduced by Freidlin and Wentzell for the study of large deviations for stationary measures of stochastic ODEs on a compact manifold, and some ideas introduced by Sowers. Some ingredients of the proof rely on rather nonstandard techniques.


\end{abstract}

\medskip
\tableofcontents

\setcounter{section}{-1}

\bigskip
\section{Introduction}
We study the large deviations for the family of probability measures $(\mu^\es)_{\es>0}$, where $\mu^\es$ stands for the invariant measure of the Markov process generated by the flow of equation
\be\label{0.1}
\p_t^2u+\gamma \p_tu-\de u+f(u)=h(x)+\sqrt{\es}\, \vartheta(t,x),\q [u(0),\dt u(0)]=[u_0, u_1].
\ee
The space variable $x$ belongs to a bounded domain $D\subset\rr^3$ with a smooth boundary, and the equation is supplemented with the Dirichlet boundary condition. The nonlinear term $f$ satisfies the dissipativity and growth conditions that are given in the next section. The force $\vartheta(t,x)$ is a colored white noise of the form
\be\label{1.63}
\vartheta(t,x)=\sum_{j=1}^\infty b_j \dt\beta_j(t)e_j(x).
\ee
Here $\{\beta_j(t)\}$ is a sequence of independent standard Brownian motions, $\{e_j\}$ is an orthonormal basis in $L^2(D)$ composed of the eigenfunctions of the Dirichlet Laplacian, and $\{b_j\}$ is a sequence of positive numbers that goes to zero sufficiently fast (see \eqref{1.57}). 
The initial point $[u_0,u_1]$ belongs to the phase space $\h=H^1_0(D)\times L^2(D)$. Finally, $h(x)$ is a function in $H^1_0(D)$ and satisfies a genericity assumption given in next section.  As it was shown in \cite{DM2014}, under the above hypotheses, the Markov process corresponding to \eqref{0.1} has a unique stationary measure $\mu^\es$ which exponentially attracts the law of any solution.

\medskip
Here we are interested in the asymptotic behavior of the family $(\mu^\es)$ as $\es$ goes to zero. We show that this family satisfies the large deviations principle (LDP), which means that there is a function that  describes precisely the logarithmic asymptotics  of $(\mu^\es)$ as the amplitude of the noise tends to zero. More formally, we have the following theorem which is part of the main result of this paper.
\begin{mt}
Let the above conditions be satisfied. Then there is a function $\vvv:\h\to [0,\iin]$ with compact level sets such that we have
\be\label{9.43}
-\inf_{\uu\in \dot\Gamma}\vvv(\uu)\le\liminf_{\es\to 0}\es\ln \mu^\es(\Gamma)\le\limsup_{\es\to 0}\es\ln \mu^\es(\Gamma)\le -\inf_{\uu\in \bar\Gamma}\vvv(\uu),
\ee
where $\Gamma$ is any Borel subset of $\h$, and we denote by $\dot\Gamma$ and $\bar\Gamma$ its interior and closure, respectively.
\end{mt}

\medskip
Before outlining the main ideas behind the proof of this result, we discuss some of the earlier works concerning the large deviations of stochastic PDEs. There is now a vast literature on this subject and the theory is developed in several directions. The most studied among them are the large deviations for the laws of trajectories of stochastic PDEs with vanishing noise. The SPDEs considered in this context include the reaction-diffusion equation \cite{sowers-1992a, CerRoc2004}, the 2D Navier-Stokes equations \cite{Chang1996, SrSu2006}, the nonlinear Schr\"odinger equation \cite{Gautier2005-2}, the Allen-Cahn equation \cite{HairWeb2014}, the quasi-geostrophic equations \cite{LiuRocZhuChan2013}, equations with general monotone drift \cite{Liu2010}, and scalar conservation laws \cite{Mariani2010}. See also the papers \cite{KalXi1996, CheMil1997, CarWeb1999, CM-2010} for results in a more abstract setting that cover a wide class of SPDEs including 2D hydrodynamical type models. Another direction is the study of exit problems for trajectories of stochastic PDEs. The results include \cite{Peszat1994, CheZhiFre2005, Gautier2008, FreKor2010, CerSal2014, BreCerFre2014}.

\medskip
The situation is completely different if we restrict our attention to the results devoted to the small-noise large deviations for stationary measures of stochastic PDEs. To the best of our knowledge, the only papers where the LDP is derived in this context are those by Sowers \cite{sowers-1992b} and Cerrai-R\"ockner \cite{CeRo2005}. These two important works are devoted to the LDP for stationary measures of the reaction-diffusion equation. In the first of them, the force is a non-Gaussian perturbation, while the second one deals with a multiplicative noise. In both papers, the origin is a unique equilibrium of the unperturbed equation and the noise is assumed to be sufficiently irregular with respect to the space variable. To the best of our knowledge, the present paper provides the first result of large deviations  for stationary measures of stochastic PDEs in the case of nontrivial limiting dynamics. Moreover, the random force $\vartheta(t,x)$ is spatially regular in our case. Both these facts create substantial additional problems which are discussed below.

\medskip
We now turn to outlining some ideas of the proof of our main result and describing the main novelty of this paper. 
Our proof relies on a development of Freidlin-Wentzell's approach.  In order to explain it, we briefly recall the original method, which relies on three main steps. The first one consists of establishing some large deviations estimates for the family of discrete-time Markov chains $(Z_n^\es)$ on the boundary. Next, one considers the family $(\lm^\es)$ of stationary measures of these chains and shows that similar estimates hold for $(\lm^\es)$. The final step is to use the Khasminskii formula to reconstruct the measure $\mu^\es$ through $\lm^\es$ and use the estimates derived for the latter in the second step, to get the LDP for $(\mu^\es)$. It turns out that in the PDE setting, this method breaks down already in the second step. Indeed, the existence of stationary measure $\lm^\es$ for the chain on the boundary is a highly nontrivial fact in this case, since on the one hand the Doob theorem cannot be applied, on the other hand this chain does not  possess the Feller property in case of a smooth random force. Moreover, even if we assume that the stationary measure exists, the classical argument does not allow to derive the LDP in this case, since the compactness of the phase space is needed.

\medskip
 To overcome these problems, we introduce a notion of \mpp{quasi-stationary measure}, which is, informally speaking, a measure that is stationary but is not supposed to be $\sigma$-additive. We show that any discrete-time Markov chain possesses such a state, thus ensuring existence of stationary measure $\lm^\es$ for the chain on the boundary in this weaker sense. It turns out that at this point (this corresponds to the second step mentioned above) the argument developed by Freidlin and Wentzell does not use the $\sigma$-additivity of $\lm^\es$, and once the necessary estimates for $(Z_n^\es)$ are obtained, they imply similar bounds for $(\lm^\es)$. Here our use of the classical technique ends, and the proof goes in a completely different direction. The reason for this is that the initial measure $\mu^\es$ cannot be reconstructed through $\lm^\es$, since, unlike the previous step, here we do need the $\sigma$-additivity of the measure $\lm^\es$. To handle this new problem,  we use the estimates obtained for $(\lm^\es)$ together with the mixing property of $\mu^\es$ established in \cite{DM2014}, to construct an auxiliary finitely additive measure $\hat\mu^\es$ defined on Borel subsets of $\h$ that satisfies
\be\label{9.65''}
\mu^\es(\dt\Gamma)\le\hat\mu^\es(\dt\Gamma)\le\hat\mu^\es(\bar\Gamma)\le\mu^\es(\bar\Gamma)\q\text{ for any }\Gamma\subset \h
\ee
and such that the family $(\hat\mu^\es)$ obeys some large deviations estimates on the balls.  The proof of the upper bound in these estimates is not a problem. The lower bound relies on an additional new ingredient, namely the notion of \mpp{stochastic stability} of a set.

\medskip
We say that a set $E\subset\h$ is \mpp{stochastically stable} if we have \,\footnote{Let us note that in the case when it is known a priori that a family $(\mu^\es)$ satisfies the LDP with a rate function $\vvv$, then a set $E$ is stochastically stable if and only if its closure has a nonempty intersection with the kernel of $\vvv$.}
$$
\lim_{\es\to 0}\es\ln\mu^{\es}(E_\eta)=0 \q\q \text{ for any }\eta>0,
$$
where $E_\eta$ stands for the open $\eta$-neighborhood of $E$ in $\h$. 

\medskip
We use it in the following context. Let us denote by $\E\subset\h$ the set of stationary flows $\uu=[u, 0]$ of the unperturbed equation
\be\label{1.51}
\p_t^2u+\gamma \p_tu-\de u+f(u)=h(x).
\ee
\begin{lemma}\label{9.47}
The set $\E$ of equilibria of \ef{1.51} is stochastically stable.
\end{lemma}
This result allows to prove the above mentioned lower bound and to complete the proof of large deviations on balls for the family $(\hat\mu^\es)$. Inequality \ef{9.65''} implies that similar result holds for the family $(\mu^\es)$ of stationary measures. The final step is to prove that this family is exponentially tight and to show that this combined with the above large deviations estimates implies the LDP.

\medskip
 We now present another essential component of the proof which allows, in particular, to get exponential tightness and also prove Lemma \ref{9.47}. Let us consider the semigroup $S(t):\h\to \h$ corresponding to \ef{1.51} and denote by $\aaa$ its global attractor.

\begin{proposition}[A priori upper bound]\label{9.18}
Under the above hypotheses, there is a function $V_\aaa:\h\to [0, \iin]$ with compact level sets and vanishing only on the attractor $\aaa$ that provides the large deviations upper bound for the family $(\mu^\es)$, that is we have
\be\label{9.19}
\limsup_{\es\to 0}\es\ln \mu^{\es}(F)\leq-\inf_{\uu\in F} V_\aaa(\uu)\q \text{ for any } F\subset\h \text{ closed}.
\ee
In particular, the family $(\mu^\es)$ is exponentially tight and any of its weak limits is concentrated on the set $\aaa$. 
\end{proposition}
Let us mention that function $V_\aaa$ has an explicit interpretation in terms of the quasipotential. Namely, for any $\uu\in \h$, $V_\aaa(\uu)$ represents the minimal energy needed to reach arbitrarily small neighborhood of $\uu$ from the global attractor in a finite time. It should be emphasized that once the main result of the paper is established, this proposition will lose its interest, since, in general, $V_\aaa$ is not the function that governs the LDP of the family $(\mu^\es)$, and that much more is proved concerning weak limits of $(\mu^\es)$. Let us mention also that some ideas of the proof of Proposition \ref{9.18} are inspired by \cite{sowers-1992b}.

\medskip
At the end of this section, let us point out that when equation \ef{1.51} has a unique equilibrium, Proposition \ref{9.18} is sufficient to derive the LDP, and in this particular case there is no need to use the Freidlin-Wentzell theory and the above scheme. Indeed, we first note that in this case the attractor $\aaa$ is a singleton $\{\hat\uu\}$, where $\hat\uu=[\hat u, 0]$ is the equilibrium position. In view of Proposition \ref{9.18}, the family $(\mu^\es)$ is tight and any weak limit of it is concentrated on $\aaa=\{\hat\uu\}$. Therefore, $\mu^\es$ weakly converges to the Dirac measure concentrated at $\hat\uu$. A simple argument (see Section \ref{9.20}) shows that this convergence and the fact that $\aaa=\{\hat\uu\}$ imply that the function $V_\aaa$ provides also the large deviations lower bound for $(\mu^\es)$. Thus, in the case of the trivial dynamics, the function $V_\aaa$ governs the LDP of the family $(\mu^\es)$. We note also that this is the only case when that happens.

\medskip
The paper is organized as follows. In Section \ref{1.76}, we state the main result and present the scheme of its proof. In Section \ref{9.45}, we establish bounds for one-step transition probabilities for the chain on the boundary. The next two sections are devoted to the proof of large deviations estimates on the balls for $(\mu^\es)$. In Section \ref{1.77}, we establish Proposition \ref{9.18}. Finally, the appendix contains some auxiliary results used in the main text.

\bigskip
{\bf Acknowledgments}. I am grateful to my supervisor Armen Shirikyan for attracting my attention to this problem, and for numerous stimulating discussions. I also thank L. Koralov and G. Raugel for discussions and some useful references. This research was carried out within the MME-DII Center of Excellence (ANR 11 LABX 0023 01) and partially supported by the ANR grant STOSYMAP (ANR 2011 BS01 015 01).

\bigskip
\section{Main result and scheme of its proof}\label{1.76}
In this section we state the main result of the paper and outline its proof. We start by recalling the notion of large deviations.
\subsection{Large deviations: equivalent formulations}
Let $\zzz$ be a Polish space. A functional $\mathfrak{I}$ defined on $\zzz$ and with range in $[0,\iin]$ is  called a (good) rate function if it has compact level sets, which means that the set $\{\mathfrak{I}\leq M\}$ is compact in $\zzz$ for any $M\geq 0$. Let $(\mathfrak{m}^\es)_{\es>0}$ be a family of probability measures on $\zzz$. The family $(\mathfrak{m}^\es)_{\es>0}$ is said to satisfy the large deviations principle in $\zzz$ with rate function $\mathfrak{I}:\zzz\to [0,\iin]$ if the following two conditions hold.

\bi
\item{\it Upper bound}
\ei
\be\label{7.1}
\limsup_{\es\to 0}\es\ln \mathfrak{m}^{\es}(F)\leq-\inf_{z\in F} \mathfrak{I}(z)\q \text{ for any } F\subset\zzz \text{ closed}.
\ee
This inequality is equivalent to the following (e.g., see Chapter 12 of \cite{DZ1992}). For any positive numbers $\De, \De'$ and $M$ there is $\es_*>0$ such that
\be\label{7.2}
\mathfrak{m}^\es(z\in\zzz: d_{\zzz}(z,\{\mathfrak{I}\leq M\})\geq\De)\leq\exp(-(M-\De')/\es)\q\text{ for } \es\leq\es_*.
\ee

\bi
\item{\it Lower bound}
\ei
\be\label{7.3}
\liminf_{\es\to 0}\es\ln\mathfrak{m}^{\es}(G)\geq-\inf_{z\in G} \mathfrak{I}(z)\q \text{ for any } G\subset\zzz \text{ open}.
\ee
This is equivalent to the following. For any $z_*\in\zzz$ and any positive numbers $\eta$ and $\eta'$ there is $\es_*>0$ such that
\be\label{7.4}
\mathfrak{m}^\es(z\in\zzz: d_\zzz(z, z_*)\le\eta)\geq\exp(-(\mathfrak{I}(z_*)+\eta')/\es)\q\text{ for } \es\leq\es_*.
\ee
The family of random variables $(\mathfrak{X}^{\es})_{\es>0}$ in $\zzz$ is said to satisfy the LDP with rate function $\mathfrak{I}$, if so does the family of their laws.

\subsection{Main result}
Before stating the main result, let us make the precise hypotheses on the nonlinearity and the coefficients entering the definition of $\vartheta(t)$. We suppose that function $f$ satisfies the growth restriction
\be\label{1.54}
|f''(u)|\leq C(|u|^{\rho-1}+1)\q u\in\rr,
\ee  
where $C$ and $\rho<2$ are positive constants, and the dissipativity conditions
\be\label{1.56}
F(u)\geq -\nu u^2-C,\q\q f(u)u- F(u)\geq-\nu u^2-C \q u\in\rr\,,
\ee
where $F$ is the primitive of $f$, $\nu\leq (\lm_1\wedge\gamma)/8$ is a positive constant, and $\lm_j$ stands for the eigenvalue corresponding to $e_j$. The coefficients $b_j$ are positive numbers satisfying
\be \label{1.57}
\BBB_1=\sum_{j=1}^\iin\lm_j b_j^2<\infty. 
\ee
Recall that we denote by $\E\subset\h$ the set of stationary flows $\uu=[u, 0]$ of equation \ef{1.51}.
It is well known that generically with respect to $h(x)$, the set $\E$ is finite (see Section \ref{9.21} for more details). We assume that $h(x)$ belongs to this generic set, so that there are finitely many equilibria, and we write $\E=\{\hat\uu_1, \ldots, \hat\uu_\ell\}$. Recall that the equilibrium $\hat \uu$ is called Lyapunov stable if for any $\eta>0$ there is $\De>0$ such that any flow of \ef{1.51} issued from the $\De$-neighborhood of $\hat\uu$ remains in the $\eta$-neighborhood of $\hat\uu$ for all time. We shall denote by $\E_s\subset \E$ the set of Lyapunov stable equilibria.

\medskip
The following theorem is the main result of this paper.
\begin{theorem}\label{1.58}
Let the above conditions be satisfied. Then the family $(\mu^\es)$ satisfies the large deviations principle in $\h$. Moreover, the corresponding rate function can vanish only on the set $\E_s\subset\{\hat\uu_1, \ldots, \hat\uu_\ell\}$ of Lyapunov stable equilibria of \eqref{1.51}. In particular, $(\mu^\es)$ is exponentially tight and any weak limit of this family is concentrated on $\E_s$. 
\end{theorem}

Let us mention that in the case when there is only one stable equilibrium $\hat\uu$ among $\{\hat \uu_1, \ldots, \hat\uu_\ell\}$ (which is the case, for example, when $\ell\leq 2$) the description of the rate function $\vvv:\h\to[0, \iin]$ that governs the LDP is quite explicit in terms of energy function (quasipotential). Namely, given $\uu$ in $\h$, $\vvv(\uu)$ represents the minimal energy needed to reach arbitrarily small neighborhood of point $\uu$ from $\hat\uu$ in a finite time. In the particular case, when the limiting equation of a stochastic PDE possesses a unique equilibrium that is globally asymptotically stable, this type of description was obtained for stochastic reaction-diffusion equation in papers \cite{sowers-1992b} and \cite{CeRo2005}.

\subsection {Scheme of the proof} 

In what follows we admit Proposition \ref{9.18}, whose proof is given in Section \ref{1.77}.

\bigskip
{\it Construction of the rate function.} We first introduce some notation following \cite{FW2012}; see Section 2 of Chapter 6 of that book. Given $\ell\in\nn$ and $i\leq\ell$, 
denote by $G_{\ell}(i)$ the set of graphs consisting of arrows 
$$
(m_1\to m_2 \to\cdots \to m_{\ell-1}\to m_\ell)
$$
such that
$$
\{m_1,\ldots, m_\ell\}=\{1,\ldots, \ell\}\q\text{ and } m_\ell=i.
$$
Further, let us introduce 
\be\label{8.47}
W_\ell(\hat\uu_i)=\min_{\mathfrak{g}\in G_{\ell}(i)}\sum_{(m\to n)\in \mathfrak{g}} V(\hat\uu_m, \hat\uu_n),
\ee
where $V(\uu_1, \uu_2)$ is the minimal energy needed to reach arbitrarily small neighborhood of $\uu_2$ from $\uu_1$ in a finite time (see \eqref{8.49} for the precise definition). The rate function $\vvv:\h\to [0, \iin]$ that governs the LDP of the family $(\mu^\es)$ is given by
\be\label{8.45}
\vvv(\uu)=\min_{i\leq \ell}[W_\ell(\hat \uu_i)+V(\hat\uu_i, \uu)]-\min_{i\leq\ell}W_\ell(\hat\uu_i).
\ee

Let us mention that when calculating these minima, we can restrict ourselves to considering only those $i$ for which $\hat\uu_i$ is Lyapunov stable.

\bigskip
{\it Markov chain on the boundary.}  What follows is a modification of a construction introduced in \cite{FW2012} (see Chapter 6) which itself is a variation of an argument used in \cite{Khas2011}. Let $\hat\uu_1,\ldots, \hat\uu_\ell$ be the stationary points of $S(t)$. Let us fix any $\uu\in\h\backslash\{\hat\uu_1,\ldots,\hat\uu_{\ell}\}$ and write $\hat\uu_{\ell+1}=\uu$. Given any $\rho_*>0$ and $0<\rho'_1<\rho_0'<\rho_1<\rho_0<\rho_*$, we use the following construction. For $i\leq \ell$, we denote by $g_i$ and $\tilde g_i$ the open $\rho_1$- and $\rho_0$-neighborhoods of $\hat\uu_i$, respectively. Similarly, we denote by $g_{\ell+1}$ and $\tilde g_{\ell+1}$, respectively, the $\rho_1'$- and $\rho_0'$-neighborhoods of $\hat\uu_{\ell+1}$. Further, we denote by $g$ and $\tilde g$ the union over $i\leq \ell+1$ of $g_i$ and $\tilde g_i$, respectively. For any $\es>0$ and $\vv\in \h$ let $S^\es(t;\vv)$ be the flow at time $t$ of \ef{0.1} issued from $\vv$. Let $\sigma_0^\es$ be the time of the first exit of the process $S^\es(t;\cdot)$ from $\tilde g$, and let $\tau_1^\es$ be the first instant after $\sigma_0^\es$ when $S^\es(t;\cdot)$ hits the boundary of $g$. Similarly, for $n\geq 1$ we denote by $\sigma_n^\es$ the first instant after $\tau_n^\es$ of  exit from $\tilde g$ and by $\tau_{n+1}^\es$ the first instant after $\sigma_n$ when $S^\es(t;\cdot)$ hits $\p g$. Let us mention that all these Markov times are almost surely finite and, moreover, have finite exponential moments (see \ef{9.22}). We consider the Markov chain on the boundary $\p g$ defined by $Z_n^\es (\cdot)= S^\es(\tau^\es_n,\cdot)$. We shall denote by $\tilde P^\es(\vv,\Gamma)$ the one-step transition probability of the chain $(Z_n^\es)$, that is 
$$
\tilde P^\es(\vv,\Gamma)=\pp(S^\es(\tau_1^\es;\vv)\in \Gamma)\q\text{ for any } \vv\in\p g\,\text{ and } \Gamma\subset\p g.
$$
The first step is a result for quasi-stationary measure $\lm^\es$ of $\tilde P_\es(\vv, \Gamma)$. We confine ourselves to announcing the result and refer the reader to Section \ref{9.82} for the definition of this concept.
\begin{proposition}\label{9.5'}
For any $\beta>0$ and $\rho_*>0$ there exist $0<\rho'_1<\rho_0'<\rho_1<\rho_0<\rho_*$ such that for all $\es<<1$ (i.e., sufficiently small), we have
\be\label{8.27'}
\exp(-(\vvv(\hat\uu_j)+\beta)/\es)\leq \lm^\es({g_j})\leq \exp(-(\vvv(\hat\uu_j)-\beta)/\es),
\ee
where inequalities hold for any $i, j\leq \ell+1, i\neq j$.
\end{proposition}
This allows to show that for $\es<<1$, there is a finitely additive measure on $\h$ satisfying \ef{9.65''} and such that
\be\label{9.81'}
\exp(-(\vvv(\hat\uu_j)+\beta)/\es)\leq \hat\mu^\es(g_j)\leq \exp(-(\vvv(\hat\uu_j)-\beta)/\es).
\ee
 As a direct corollary of these relations, we get the following result. 
\begin{proposition}\label{9.4} For any $\beta>0$ and $\rho_*>0$ there exist $0<\rho'_1<\rho'_0<\rho_1<\rho_0<\rho_*$ such that for any $j\leq \ell+1$ and $\es<<1$, we have
\begin{align}
\mu^\es(g_j)&\leq \exp(-(\vvv(\hat\uu_j)-\beta)/\es)\label{9.79},\\
\mu^\es(\bar g_j)&\ge \exp(-(\vvv(\hat\uu_j)+\beta)/\es)\label{9.80}.
\end{align}
\end{proposition}
The passage from Proposition \ref{9.5'} to \ref{9.4} is the most involved part of the paper and construction of $\hat\mu^\es$ is the main idea behind its proof. Without going into details, we describe in few words another key ingredient of the proof, namely Lemma \ref{9.47}.
\begin{definition}
We shall say that a set $E\subset\h$ is \mpp{stochastically stable} or \mpp{stable with respect to}  $(\mu^\es)$ if we have 

$$
\lim_{\es\to 0}\es\ln\mu^{\es}(E_\eta)=0 \q\q \text{ for any }\eta>0,
$$
where $E_\eta$ stands for the open $\eta$-neighborhood of $E$ in $\h$. If the above relation holds only along some sequence $\es_j\to 0$ (that is $\es$ replaced by $\es_j$), we shall say that $E$ is \mpp{stable with respect to} $(\mu^{\es_j})$.
\end{definition}
Let us show how to derive Lemma \ref{9.47} using Proposition \ref{9.18}. We then show how the same proposition combined with \eqref{9.79}-\ef{9.80} implies the LDP.
 We admit the following result established in the appendix.
\begin{lemma}\label{1.201}
Let $\ooo$ be a heteroclinic orbit of $S(t)$ and let $\uu_1\in\ooo$. Suppose that $\uu_1$ is stable with respect to $(\mu^{\es_j})$ for some sequence $\es_j\to 0$. Then so is any point $\uu_2$ lying on $\ooo$ after $\uu_1$ (in the direction of the orbit).
\end{lemma}
Let us mention that we consider the endpoints of an orbit as its elements, and when saying $\uu$ is stable with respect to $(\mu^{\es_j})$ we mean that so is the set $\{\uu\}$. Now let us assume that Lemma \ref{9.47} is not true. Then we can find two positive constants $a$ and $\eta$, and a sequence $\es_j$ going to zero such that
\be\label{8.20}
\mu^{\es_j}(\E_\eta)=\sum_{j=1}^\ell \mu^{\es_j}(B(\hat\uu_j,\eta)) \leq\exp(-a/\es_j)\q\text{ for all }j\geq 1,
\ee
where $B(\uu, r)$ stands for the open ball in $\h$ of radius $r$ and centered at $\uu$. 
By Proposition \ref{9.18}, the sequence $(\mu^{\es_j})$ is tight and any weak limit of it is concentrated on $\aaa$. So, up to extracting a subsequence, we can assume that $\mu^{\es_j}\rightharpoonup \mu_*$, and $\mu_*$ is concentrated on $\aaa$. By Theorem \ref{Th-attractor}, the global attractor $\aaa$ consists of points $(\hat\uu_i)_{i=1}^n$ and joining them heteroclinic orbits. Let $\uu_*$ be a point lying on such an orbit that belongs to the support of $\mu_*$. By the portmanteau theorem, we have
$$
\liminf_{j\to \iin}\mu^{\es_j}(B(\uu_*,r))\ge\mu_*(B(\uu_*,r))>0\q\q\text{ for any }r>0.
$$
Therefore, the point $\uu_*$ is stable with respect to $(\mu^{\es_j})$. On the other hand, it follows from the previous lemma that so are all points of the attractor that lie on that orbit  after $\uu_*$. In particular so is the endpoint of $\ooo$, which is in contradiction with \eqref{8.20}. The proof of Lemma \ref{9.47} is complete.

\bigskip
{\it Derivation of the LDP.} 
We claim that the hypotheses of Lemma \ref{9.3} are satisfied for the family $(\mu^\es)_{\es>0}$ and rate function $\vvv$. Indeed, let $\beta$ and $\rho_*$ be two positive constants and let $\uu$ be any point in $\h$. If $\uu$ is not a stationary point, we denote $\hat\uu_{\ell+1}=\uu$ and use Proposition \ref{9.4} to find $\rho_1'<\rho_1<\rho_*$ such that we have \eqref{9.79}-\ef{9.80} and we set $\tilde\rho(\uu)=\rho_1'$. Otherwise ($\uu$ is stationary), we take any non stationary point $\uu'$ and denote $\hat\uu_{\ell+1}=\uu'$. We once again use Proposition \ref{9.4} to find $\rho_1'<\rho_1<\rho_*$ such that we have \eqref{9.79}-\ef{9.80} and we set $\tilde\rho(\uu)=\rho_1$. Let us note in this case ($\uu$ is stationary) the choice of $\uu'$ is not important due to the fact that we are interested in the asymptotic behavior of $(\mu^\es)$  only in the neighborhood of $\uu$, and we add a new point $\hat\uu_{\ell+1}=\uu'$ only to be consistent with Proposition \ref{9.4}. Thus, the hypotheses of Lemma \ref{9.3} are satisfied and the family $(\mu^\es)$ satisfies the LDP in $\h$ with rate function $\vvv$.

\bigskip
\section{Proof of Theorem \ref{1.58}}\label{9.45}
The present section is devoted to the proof of the main result of this paper. We admit Proposition \ref{9.18}, which is proved in the next section, and following the scheme presented above establish Theorem \ref{1.58}. We shall always assume that the hypotheses of this theorem are satisfied. 
 
\subsection{Construction of the rate function}\label{1.64}
Here we define the function $V:\h\times\h\to [0, \iin]$ entering relation \ef{8.45} and function $V_\aaa:\h\to [0, \iin]$ from Proposition \ref{9.18}. We first introduce some notation. 
For any $t\geq 0$, $\vv\in\h$ and $\ph\in L^2(0,T;L^2(D))$, let us denote by $S^\ph(t;\vv)$ the flow at time $t$ of equation 
\be\label{7.10}
\p_t^2 u+\gamma \p_t u-\de u+f(u)=h(x)+\ph(t,x)
\ee
issued from $\vv$.
Let $H_{\vartheta}$ be the Hilbert space defined by \be\label{1.34}
H_{\vartheta}=\{v\in L^2(D): |v|_{H_\vartheta}^2=\sum_{j=1}^\iin b_j^{-2}\,(v,e_j)^2<\iin\}.
\ee
For a trajectory $\uu_\cdot\in C(0,T;\h)$ we introduce
\be\label{8.50}
I_T(\uu_\cdot)=J_T(\ph):=\f{1}{2}\int_0^T |\ph(s)|_{H_\vartheta}^2\dd s,
\ee
if there is $\ph\in L^2(0,T;H_\vartheta)$ such that $\uu_\cdot=S^\ph(\cdot;\uu(0))$, and $I_T(\uu_\cdot)=\iin$ otherwise.
We now define $V:\h\times\h\to [0, \iin]$ by
\be\label{8.49}
V(\uu_1,\uu_2)=\lim_{\eta\to 0}\inf\{I_T(\uu_\cdot); T>0, \uu_\cdot\in C(0,T;\h): \uu(0)=\uu_1, \uu(T)\in B(\uu_2,\eta) \}.
\ee

\medskip
Let us note that this limit (finite or infinite) exists, since the expression written after the limit sign is monotone in $\eta>0$. As we mentioned in previous section, $V(\uu_1, \uu_2)$ represents the minimal energy needed to reach arbitrarily small neighborhood of $\uu_2$ from $\uu_1$ in a finite time.

\bigskip{\it Remark.} The definition of the quasipotential $V$ using this filtration in $\eta$ rather than  taking directly $\eta=0$ is explained by the lack of the exact controllability of the NLW equation by a regular force, and \ef{8.49} ensures the lower semicontinuity of function $\vvv$ given by \ef{8.45}.

\medskip
The function $V_\aaa:\h\to [0, \iin]$ entering Proposition \ref{9.18} is defined by
\be\label{9.6}
V_\aaa(\uu_*)=\inf_{\uu_1\in \aaa} V(\uu_1, \uu_*).
\ee

\medskip
 Notice that the compactness of level sets of $V_\aaa$ implies that $\vvv$ has relatively compact level sets. Combining this with lower semicontinuity of $\vvv$ (the proof of this fact is identical to that of $V_\aaa$, see Section \ref{1.77}), we see that $\vvv$ is a rate function in $\h$.
In what follows, the space $\h$ will be endowed with the norm
\be\label{9.11}
|\uu|_\h^2=\|\g u_1\|^2+\|u_2+\al u_1\|^2 \q\text{ for }\uu=[u_1, u_2]\in\h, \,\footnote{\, We denote by $\|\cdot\|_s$ the $H^s$-norm of a vector, and $\|\cdot\|$ stands for the $L^2$-norm.}
\ee
where $\al>0$ is a small parameter.

\subsection{Markov chain on the boundary}
In this section we establish a result that implies Proposition \ref{9.5'}. For the proof of this implication, see Chapter 6 of \cite{FW2012}; the only difference here is that $\lm^\es$ is not necessarily $\sigma$-additive, which does not affect the proof. 

\medskip
Recall that we denote by $V(\hat\uu_i, \hat\uu_j)$ the minimal energy needed to reach any neighborhood of $\hat\uu_j$ from $\hat\uu_i$ in a finite time. In what follows, we shall denote by $\tilde V(\hat\uu_i, \hat\uu_j)$ the energy needed to reach any neighborhood $\hat\uu_j$ from $\hat\uu_i$ in a finite time without intersecting any $\hat\uu_k$, for $k\le\ell+1$ different from $i$ and $j$.
\begin{proposition}\label{9.5}
For any positive constants $\beta$ and $\rho_*$ there exist $0<\rho'_1<\rho_0'<\rho_1<\rho_0<\rho_*$ such that for all $\es<<1$, we have
\be\label{8.27}
\exp(-(\tilde V(\hat\uu_i,\hat\uu_j)+\beta)/\es)\leq \tilde P^\es(\vv,\p g_j)\leq \exp(-(\tilde V(\hat\uu_i,\hat\uu_j)-\beta)/\es),
\ee
where inequalities hold for any $\vv\in \p g_i$ and any $i, j\leq \ell+1, i\neq j$.
\end{proposition}

\medskip
{\it Comment.} In what follows, when proving this type of inequalities, we shall sometimes derive them with $\beta$ replaced by $C\beta$, where $C\ge 1$ is an absolute constant. Since $\beta>0$ can be taken arbitrarily small, these bounds are equivalent and we shall use this without further stipulation.

\medskip
{\it Derivation of the lower bound.} We assume that $\tilde V(\hat\uu_i,\hat\uu_j)<\iin$, since otherwise there is nothing to prove. We shall first establish the bound for $i\leq \ell$.
We need the following result, whose proof is given at the end of this section.
\begin{lemma}\label{8.28}
There exists $\tilde\rho>0$ such that for any $0<\rho_2<\rho_1<\tilde\rho$ we can find a finite time $T>0$ depending only on $\rho_1$ and $\rho_2$ such that for any point $\vv\in \bar B(\hat\uu_i,\rho_1)$, $i\leq \ell$, there is an action $\ph_{\vv}$ defined on the interval $[0, T]$ with energy not greater than $\beta$ such that we have
\be\label{8.25}
S^{\ph_\vv}(t;\vv)\in \bar B(\hat\uu_i,\rho_1)\q\text{ for } t\in[0,T] \q\text{ and }\q S^{\ph_\vv}(T;\vv)\in \bar B(\hat\uu_i,\rho_2/2).
\ee
\end{lemma}
By definition of $\tilde V$, for $\rho_*>0$ sufficiently small and $\rho_1'<\rho_*$, we can find a finite time $\tilde T>0$ and an action $\tilde \ph$ defined on $[0, \tilde T]$ with energy smaller than $\tilde V(\hat\uu_i,\hat\uu_j)+\beta$ such that $S^{\tilde\ph}(\tilde T;\hat\uu_i)\in B(\hat\uu_j, \rho_1'/4)$ and the curve $S^{\tilde\ph}(\tilde T;\hat\uu_i)$ does not intersect $\rho_*$-neighborhood of $\hat\uu_k$ for $k\neq i, j$ (note that if a trajectory does not intersect $\hat\uu_k$ then it also does not intersect some small neighborhood of $\hat\uu_k$). Since $\rho_*>0$ can be taken arbitrarily small, we may assume that $\rho_*\le\tilde\rho$, where $\tilde\rho$ is the constant from the above lemma. Let $\rho_2<\rho_*$ be so small that for any $\vv\in \bar B(\hat\uu_i, \rho_2)$ we have $S^{\tilde\ph}(\tilde T;\vv)\in B(\hat\uu_j,\rho'_1/2)$. 
We take any $\rho_1\in(\rho_2, \rho_*)$ and use the following construction. For any $\vv\in \bar g_i$, we denote by $\tilde\ph_{\vv}$ the action defined on $[0, T+\tilde T]$ that coincides with $\ph_\vv$ on $[0, T]$ and with $\tilde\ph$ on $[T, T+\tilde T]$. Let us note that for any $\vv\in\bar g_i$, we have
\be\label{8.29}
I_{T+\tilde T}(S^{\tilde\ph_\vv}(\cdot;\vv))=J_{T+\tilde T}(\tilde\ph_\vv)=J_T(\ph_\vv)+J_{\tilde T}(\tilde\ph)\leq \tilde V(\hat\uu_i,\hat\uu_j)+2\beta.
\ee
Now let us take any $\rho_0\in (\rho_1,\rho_*)$ and denote by $\De$ any positive number that is smaller than $\min\{\rho_0-\rho_1, \rho_2/2, \rho_1'/2\}$. Then we have the following: if the trajectory $S^\es(t;\vv)$ is in the $\De$-neighborhood of $S^{\tilde\ph_\vv}(t;\vv)$ in $C(0, T+\tilde T;\h)$ distance, then $\tau_1^\es(\vv)\leq T+\tilde T$ and $S^\es(\tau_1^\es;\vv)\in \p g_j$. Therefore, we have
$$
\inf_{\vv\in\p g_i}\tilde P^\es(\vv,\p g_j)\geq \pp(A),
$$
where we set 
$$
A=\{\om\in\Omega:\sup_{\vv\in \bar g_i} d_{C(0,T+\tilde T;\h)}(S^\es(\cdot;\vv), S^{\tilde\ph_\vv}(\cdot;\vv))<\De\}.
$$
Combining this with inequality \eqref{8.29} and Theorem \ref{6.1}, we derive the lower bound of \eqref{8.27} in the case $i\leq \ell$.

\medskip
We now show that if $\rho_0'<\rho_1$ is sufficiently small, then the lower bound is also true for $i=\ell+1$. Indeed, let $\tilde V(\hat\uu_{\ell+1},\hat\uu_j)<\iin$ and let $T>0$ and $\ph$ be such that $S^{\ph}(T;\hat\uu_{\ell+1})\in B(\hat\uu_j,\rho_1/4)$ and 
\be\label{8.30}
J_{T}(\ph)\leq \tilde V(\hat\uu_{\ell+1},\hat\uu_j)+\beta.
\ee
We assume that $\rho_0'$ is so small that $S^{\ph}(T;\vv)\in B(\hat\uu_j,\rho_1/2)$ for any $\vv\in \tilde g_{\ell+1}$. Let us take any $\De<\min\{\rho_0'-\rho_1', \rho_1/2\}$. Then for any $\vv\in \bar g_{\ell+1}$ if the trajectory $S^\es(t;\vv)$ is in the $\De$-neighborhood of $S^{\ph}(t;\vv)$ in $C(0, T;\h)$ distance, then we have $\tau_1^\es(\vv)\leq T$ and $S^\es(\tau_1^\es;\vv)\in \p g_j$. Therefore
$$
\inf_{\vv\in\p g_{\ell+1}}\tilde P^\es(\vv,\p g_j)\geq \pp(A'),
$$
where we set 
$$
A'=\{\om\in\Omega:\sup_{\vv\in \bar g_{\ell+1}} d_{C(0,T;\h)}(S^\es(\cdot;\vv), S^{\ph}(\cdot;\vv))<\De\}.
$$
Combining this with inequality \eqref{8.30} and Theorem \ref{6.1}, we derive the lower bound in the case $i=\ell+1$. 

\bigskip
{\it Proof of the upper bound.} We assume that $\rho_*>0$ is so small that the energy needed to move the point from $\rho_*$-neighborhood of $\hat\uu_i$ to $\rho_*$-neighborhood of $\hat\uu_j$ without intersecting any other $\hat\uu_k$ is no less than $\tilde V(\hat\uu_i,\hat\uu_j)-\beta.$ Let us denote by $\tau_g^\es$ the first instant when the process $S^\es(t,\cdot)$ hits the set $\bar g$.  Then, by the strong Markov property, we have 
\be\label{8.33}
\sup_{\vv\in\p g_i}\pp(S^\es(\tau_1^\es;\vv)\in \p g_j)\leq \sup_{\vv\in\p \tilde g_i}\pp(S^\es(\tau_g^\es;\vv)\in \p g_j).
\ee
In what follows we shall denote by $g'$ the set $g\backslash g_{\ell+1}$, i.e. the union over $i\leq\ell$ of $\rho_1$-neighborhoods of $\hat\uu_i$. We need the following result.
\begin{lemma}\label{8.31}
For any positive constants $\rho_1, R$ and $M$ there is $T>0$ such that
\be\label{8.38}
\sup_{\vv\in B_R}\pp(\tau^\es_{g'}(\vv)\geq T)\leq\exp(-M/\es)\q\text{ for }\es<<1,
\ee
where $\tau^\es_{g'}$ stands for the first hitting time of the set $\bar g'$, and $B_R$ is the closed ball in $\h$ of radius $R$ centered at the origin. 
\end{lemma}
Note that for any $\vv\in\p\tilde g_i$, we have
$$
\pp(S^\es(\tau_g^\es;\vv)\in \p g_j)\leq\pp(S^\es(\tau_g^\es;\vv)\in \p g_j, \tau_g^\es(\vv)< T)+\pp(\tau^\es_{g'}(\vv)\geq T),
$$
where we used the fact that $\tau^\es_{g}\leq\tau^\es_{g'}$. Now notice that the event under the probability sign of the first term of this sum means that the trajectory $S^\es(t;\cdot)$ issued from $\p\tilde g_i$ hits the set $\bar g_j$ over time $T$ and does not intersect any $\hat\uu_k$ for $k$ different from $i$ and $j$. It follows from Theorem \ref{6.1} that this event has a probability no greater than $\exp(-(\tilde V(\hat\uu_i,\hat\uu_j)-2\beta)/\es)$. Combining this with Lemma \ref{8.31} and inequality \eqref{8.33}, we infer
$$
\sup_{\vv\in\p g_i}\tilde P^\es(\vv,\p g_j)\leq \exp(-(\tilde V(\hat\uu_i,\hat\uu_j)-2\beta)/\es)+\exp(-M/\es)\q\text{ for }\es<<1.
$$
Since $M>0$ can be chosen arbitrarily large, we derive the upper bound.

\bigskip
{\it Proof of Lemma \ref{8.28}.}
For any $\vv\in \bar B(\hat\uu_i,\rho_1)$, let $\tilde\vv(t;\vv)$ be the flow issued from $\vv$ corresponding to the solution of
$$
\p_t^2 \tilde v+\gamma\p_t \tilde v-\de \tilde v+f(\tilde v)=h(x)+P_N[f(\tilde v)-f(\hat u_i)],
$$
where $\hat\uu_i=[\hat u_i, 0]$ and $P_N$ stands for the orthogonal projection from $L^2$ to its $N$ dimensional subspace spanned by vectors $e_1,\ldots, e_N$. Let us define $\ph_{\vv}=P_N[f(\tilde v)-f(\hat u_i)]$. Then we have $\tilde\vv(t;\vv)=S^{\ph_\vv}(t;\vv)$. Moreover, it follows from Proposition \ref{1.88} that for $N\geq N_*(|\hat\uu_i|_\h)$ we have
\be\label{8.24}
|S^{\ph_\vv}(t;\vv)-\hat\uu_i|_\h^2\leq e^{-\al t}|\vv-\hat\uu_i|_\h^2,
\ee
where $\al>0$ is the constant entering \ef{9.11}. In particular, if we take $T=2(\ln\rho_1-\ln\rho_2)/\al$ then we get \eqref{8.25}. Moreover, we have
\begin{align}
J_{T}(\ph_\vv)&=\f{1}{2}\int_0^{T}|P_N[f(\tilde v)-f(\hat u_i)]|_{H_\vartheta}^2\dd s\leq C(N)\int_0^{T}|f(\tilde v)-f(\hat u_i)|_{L^1}^2\dd s\notag\\
&\leq C_1\, C(N)\int_0^{T}(\|\tilde v\|_1^2+\|\hat u_i\|_1^2+1)\|\tilde v-\hat u_i\|_1^2\dd s\notag\\
&\leq C_2\,C(N) \int_0^{T}|\vv-\hat\uu_i|^2 e^{-\al s}\dd s\leq C_3\, C(N)\, \tilde\rho^2\leq \beta\label{5.14},
\end{align}
provided $\tilde\rho>0$ is sufficiently small.

\bigskip
{\it Proof of Lemma \ref{8.31}. Step~1: Reduction}. Let the positive constants $\rho_1, R$ and $M$ be fixed. We claim that it is sufficient to prove that for any $R'> R$ we can find positive constants $T_*$ and $a$ such that 
\be\label{8.34}
\sup_{\vv\in B_{R'}}\pp(\tau^\es_{g'}(\vv)\geq T_*)\leq\exp(-a/\es)\q\text{ for }\es<<1.
\ee
Indeed, taking this inequality for granted, let us derive \eqref{8.38}.
To this end, let us use Proposition \ref{1.96} to find $R'> R$ so large that
\be\label{8.35}
\sup_{t\geq 0}\sup_{\vv\in B_R}\pp(S^\es(t;\vv)\notin B_{R'})\leq\exp(-(M+1)/\es)\q\text{ for }\es<1.
\ee
Once such $R'$ is fixed we find $T_*>0$ and $a>0$ such that we have \eqref{8.34}.
Let us take $n\geq 1$ so large that $an>(M+1)$.
For any $k\leq n$ we introduce
$$
p_k=\sup_{\vv\in B_R}\pp(\tau^\es_{g'}(\vv)\geq k\,T_*).
$$
Then, we have
\begin{align}\label{8.37}
p_n&\leq\sup_{\vv\in B_R}\pp(\tau^\es_{g'}(\vv)\geq n\,T_*, S^\es((n-1)T_*;\vv)\in B_{R'})\notag\\
&\q+\sup_{\vv\in B_R}\pp(S^\es((n-1)T_*;\vv)\notin B_{R'})\leq q_n+\exp(-(M+1)/\es),
\end{align}
where we denote by $q_n$ the first term of this sum and we used inequality \eqref{8.35} to estimate the second one. Now note that by the Markov property, we have
\begin{align*}
q_n&=\sup_{\vv\in B_R}\e_\vv[\e(\ch_{\tau^\es_{g'}\geq n\,T_*}\ch_{S^\es((n-1)T_*;\vv)\in B_{R'}}|\fff^\es_{(n-1) T_*})]\\
&=\sup_{\vv\in B_R}\e_\vv[\ch_{\tau^\es_{g'}\geq (n-1)\,T_*} \ch_{S^\es((n-1)T_*;\vv)\in B_{R'}}\pp(\tau^\es_{g'}(S^\es((n-1)T_*;\vv))\geq T_*)]\\
&\leq \sup_{\tilde\vv\in B_{R'}}\pp(\tau^\es_{g'}(\tilde\vv)\geq T_*)\,p_{n-1}\leq\exp(-a/\es)\,p_{n-1},
\end{align*}
where $\fff^\es_t$ is the filtration corresponding to $S^\es(t;\cdot)$, and we used inequality \eqref{8.34}. Combining this with \eqref{8.37}, we derive
$$
p_n\leq \exp(-a/\es)\,p_{n-1}+\exp(-(M+1)/\es).
$$
Iterating this inequality, we infer 
$$
p_n\leq \exp(-an/\es)+(1-\exp(-a/\es))^{-1}\exp(-(M+1)/\es)\leq\exp(-M/\es).
$$
It follows that inequality \eqref{8.38} holds with $T=n\, T_*$.

\bigskip
{\it Step~2: Derivation of \eqref{8.34}}. We first show that for any positive constants $\tilde R$ and $\eta$, we have
\be\label{8.36}
\sup_{\vv\in B_{\tilde R}}l(\vv)<\iin,
\ee
where $\l(\vv)$ stands for the first instant when the deterministic flow $S(t)\vv$ hits the set $\bar\ooo(\eta)$ and where $\ooo(\eta)$ is the union over $i\leq\ell$ of $\eta$-neighborhoods of $\hat\uu_i$. Indeed, let us suppose that this is not true, and let us find $\tilde R>0$ and $\eta>0$ for which this inequality fails. 
Then, there exists a sequence $(\vv_j)\subset B_{\tilde R}$ such that 
\be\label{8.39}
l(\vv_j)\geq 2j.
\ee
For each $j\geq 1$, let us split the flow $S(t)\vv_j$ to the sum $\uu^1_j(t)+\uu^2_j(t)$, where $\uu^1_j(t)$ stands for the flow issued from $\vv_j$ of equation \eqref{1.51} with $f=h=0$. Then, for all $t\geq 0$, we have
\be\label{8.40}
|\uu^1_j(t)|_\h^2\leq e^{-\al t}|\vv_j|_\h^2, \q\q |\uu^2_j(t)|_{H^{s+1}\times H^s}\leq C_s(\tilde R),
\ee
where $s<1-\rho/2$ is any constant (e.g., see \cite{BV1992, Har85}). Using second of these inequalities, let us find $(j_k)\subset\nn$ such that the sequence $\uu^2_{j_k}(j_k)$ converges in $\h$ and denote by $\tilde \uu$ its limit. Then, in view of first inequality of \eqref{8.40}, we have 
\be\label{8.41}
S(j_k)\vv_{j_k}\to\tilde \uu \q\text{ in } \h \q\text{ as } k\to\iin.
\ee
Now let us find $t_*>0$ so large that
\be\label{8.42}
S(t_*)\tilde\uu\in \bar\ooo(\eta/2).
\ee
Note that thanks to \eqref{8.41} and continuity of $S(t)$, we have
\be\label{8.43}
S(j_k+t_*)\vv_{j_k}\to S(t_*)\tilde\uu.
\ee
Further, notice that by \eqref{8.39} we have $S(j_k+t_*)\vv_{j_k}\notin \bar\ooo(\eta)$ for $k\geq 1$ sufficiently large. This is clearly in contradiction with \eqref{8.42}-\eqref{8.43}. Inequality \eqref{8.36} is thus established.

\medskip
We are now ready to prove \eqref{8.34}. Indeed, let us use inequality \eqref{8.36} with $\tilde R=R'+1$ and $\eta=\rho_1/2$, and let us set 
\be\label{8.44}
T_*=\sup_{\vv\in B_{R'+1}} l(\vv)+1.
\ee
Let us consider the trajectories $\uu_\cdot\in C(0,T_*;\h)$ issued from $B_{R'+1}\backslash \ooo(\rho_1/2)$ and assuming their values outside $\ooo(\rho_1/2)$. Note that the family $\elll$ of such trajectories is closed in $C(0, T_*;\h)$. Therefore, the infimum 
$$
a'=\inf_{\uu_\cdot\in\elll}I_{T_*}(\uu_\cdot)
$$
is attained and is positive, since in view of \eqref{8.44} there is no deterministic trajectory $S(t)\vv$ under consideration. Now using Theorem \ref{6.1}, we see that \eqref{8.34} holds with $a=a'/2$.
The proof of Lemma \ref{8.31} is complete.

\section{Quasi-stationarity and auxiliary measure}\label{9.82}
In this section we introduce a notion of quasi-stationary measure and show that any discrete-time Markov chain possesses such a state. This will be used to construct a finitely additive measure $\hat\mu^\es$ satisfying relation \ef{9.65''} and such that the family $(\hat\mu^\es)$ satisfies some large deviations estimates (see \ef{9.78}). This in turn will imply Proposition \ref{9.4}.
\subsection{Quasi-stationary measure}
 Let $X$ be a metric space and let $b(X)$ the space of bounded Borel measurable functions on $X$ equipped with the topology of uniform convergence. We shall denote by $b^*(X)$ the dual of $b(X)$ \,\footnote{\, $b^*(X)$ can be identified  with the space $ba(X)$ of finite, finitely additive signed measures on $X$ (e.g., see Theorem IV.5.1 in \cite{Dunford-Schwartz1}).} . A linear continuous map $\pP$ from $b(X)$ into itself is called a Markov operator on $X$, if $\pP\psi\ge 0$ for any $\psi\ge 0$ and $\pP1\equiv 1$. Let $\pP^*:b^*(X)\to b^*(X)$ be the dual of $\pP$, that is  
 $$
\pP^*\lm(\psi)=\lm(\pP\psi)
$$
for any $\lm\in b^*(X)$ and $\psi\in b(X)$. We shall say that $\lm$ is \mpp{a quasi-stationary state} (or \mpp {measure)}  for $\pP$ if it satisfies the following properties: $\lm(\psi)\ge 0$ for any $\psi\ge 0$, $\lm(1)=1$, and $\pP^*\lm=\lm$, that is we have
 \be\label{9.57}
\lm(\pP\psi)=\lm(\psi)\q\text{ for any }\psi\in b(X).
\ee
To any such $\lm$ we associate a finitely additive measure defined on Borel subsets of $X$ by $\lm(\Gamma)=\lm(\ch_\Gamma)$ for $\Gamma\subset X$.
\begin{lemma}\label{lem-stationary}
Any Markov operator possesses a quasi-stationary measure.
\end{lemma}
\bp
Let $\pP$ be a Markov operator defined on a space $X$. Consider the space
$$
\mathfrak{F}=\{\lm\in b^*(X): \lm(1)= 1 \text{ and } \lm(\psi)\ge 0 \text{ for }\psi\ge 0\}
$$
endowed with weak* topology. Note that if $\lm\in \mathfrak{F}$ then its norm is equal to 1. In view of the Banach-Alaoglu theorem, $\mathfrak{F}$ is relatively compact. Moreover, it is easy to see that $\mathfrak{F}$ is also closed and convex. Since $\pP$ is a Markov operator, its dual $\pP^*$ maps $\mathfrak{F}$ into itself. Thanks to the Leray-Schauder theorem (e.g., see Chapter 14 in \cite{taylor1996}), $\pP^*$ has a fixed point $\lm\in \mathfrak{F}$, which means that $\lm$ is quasi-stationary for $\pP$.  It should be emphasized that $\lm$ is not stationary in the classical sense, since it is not necessarily $\sigma$-additive. 
\ep

\bigskip
In what follows, given $\es>0$, we shall denote by $\lm=\lm^\es$ any of quasi-stationary states of the operator $\pP=\pP^\es:b(\p g)\to b (\p g)$ defined by
\be\label{9.75}
\pP\psi(\vv)=\int_{\p g}\psi(\Zz)\tilde P_1(\vv, \Dd\Zz)\equiv \e\psi(S^\es(\tau_1^\es;\vv)).
\ee
So we have
\be\label{9.74}
\lm(\pP\psi)=\lm(\psi)\q\text{ for any }\psi\in b(\p g)
\ee
and $\lm\in b^*(\p g)$ satisfies $\lm(\psi)\ge 0$ for $\psi\ge 0$, $\lm (1)=1$. We shall always assume that $\es>0$ is so small that we have \ef{8.27}.

\subsection{Khasminskii type relation}
For each $\es>0$, let us define a continuous map $\tilde\mu=\tilde\mu^\es$ from $b(\h)$ to $\rr$ by 
\be\label{9.59}
\tilde\mu(\psi)=\lm(\elll\psi),
\ee
where $\lm=\lm^\es$ is given by \ef{9.75}-\ef{9.74}, and $\elll=\elll^\es:b(\h)\to b(\h)$ is defined by
\be\label{9.60}
\elll\psi(\vv)=\e\int_0^{\tau_1^\es}\psi(S^\es(t;\vv))\dd t.
\ee
\begin{lemma} For any $\psi\in b(\h)$, we have
\end{lemma}
\be\label{9.61}
\tilde\mu(P_s\psi)=\tilde\mu(\psi)\q\text{ for any } s\ge 0,
\ee
where $P_s=P_s^\es: b(\h)\to b(\h)$ stands for the Markov operator of the process $S^\es(\cdot)$.
\bp
We use the classical argument (see Chapter 4 in \cite{Khas2011}). Let us fix $s\ge 0$. By the Markov property, for any $\vv\in\h$, we have
$$
\e\int_0^{\tau_1^\es}\psi(S^\es(t+s;\vv))\,dt=\e\int_0^{\tau_1^\es}P_s\psi(S^\es(t;\vv))\dd t.
$$
It follows that
\begin{align*}
\tilde\mu(P_s\psi)&=\lm(\elll P_s\psi)=\lm(\e\int_0^{\tau_1^\es}P_s\psi(S^\es(t;\cdot))\dd t)=\lm(\e\int_0^{\tau_1^\es}\psi(S^\es(t+s;\cdot))\dd t)\notag\\
&=\lm(\e\int_s^{s+\tau_1^\es}\psi(S^\es(t;\cdot))\dd t)\notag\\
&=\lm(\e\int_{\tau_1^\es}^{s+\tau_1^\es}\psi(S^\es(t;\cdot))\dd t)-\lm(\e\int_0^{s}\psi(S^\es(t;\cdot))\dd t)+\tilde\mu(\psi).
\end{align*}
Conditioning with respect to $\fff_{\tau_1^\es}$ and using the strong Markov property together with \ef{9.75}-\ef{9.74}, we see that the first two terms are equal. Since $s\ge 0$ was arbitrary, we arrive at \ef{9.61}.
\ep

\subsection{Auxiliary finitely additive measure}
Let us denote by $b_0(\h)$ the space $b(\h)$ endowed with topology of uniform convergence on bounded sets in $\h$. That is, given a point $\psi\in b(\h)$ and a sequence $(\psi_n)\subset b(\h)$, we shall say that $\psi_n$ converges to $\psi$ in $b_0(\h)$  as $n\to \iin$, if for any $a>0$ we have
$$
\sup_{\uu\in B_a}|\psi_n(\uu)-\psi(\uu)|\to 0\q\text{ as }n\to \iin.
$$ 
\begin{lemma}\label{9.58}
The map $\tilde\mu$ given by \ef{9.59}-\ef{9.60} is continuous from $b_0(\h)$ to $\rr$.
\end{lemma}
\bp
For the simplicity, we shall write $\tau_1=\tau_1^\es$ and $\uu(t;\vv)=S^\es(t;\vv)$. We need to show that $\tilde\mu(\psi_n)\to 0$ for any $\psi_n$ converging to zero in $b_0(\h)$.
Since $\lm$ is continuous from $b(\p g)$ to $\rr$, it is sufficient to show that $\elll\psi_n$ goes to zero uniformly in $B_R$, where $R>0$ is so large that $g\subset B_R$. Let us fix any $\eta>0$.
Clearly, we may assume that $|\psi_n(\uu)|\le 1$ for any $n\ge 1$ and $\uu\in \h$. 
It follows from the Cauchy-Schwarz inequality and \ef{9.22} that for $R_1>0$ sufficiently large, we have 
\be\label{9.62}
\e_\vv(\tau_1\cdot\ch_{\tau_1\ge R_1})\le \eta\q\text{ for any }\vv\in B_R.
\ee
Once $R_1$ is fixed, let us use Proposition 3.2 from \cite{DM2014} to find $R_2>R_1$ such that 
\be\label{9.63}
\pp\left(\sup_{t\in [0, R_1]}|\uu(t;\vv)|\ge R_2\right)\le\eta/{R_1}\q\text{ for any }\vv\in B_R.
\ee
Now we have
\begin{align}
|\elll\psi_n(\vv)|\le\e\int_0^{\tau_1}|\psi_n(\uu(t;\vv))|\dd t &\leq \e(\ch_{\tau_1\ge R_1}\int_0^{\tau_1}|\psi_n(\uu(t;\vv))|\dd t)\notag\\
&\q +\e\int_0^{R_1}|\psi_n(\uu(t;\vv))|\dd t:=I_1+I_2.\label{9.64}
\end{align}
Let us note that in view of \ef{9.62}, we have $I_1\le\eta$. Further, since $\psi_n$ converges to zero in $b_0(\h)$, we can find $n_*(\eta)\ge 1$ such that for all $n\ge n_*(\eta)$, we have
$$
\sup_{\uu\in B_{R_2}}|\psi_n(\uu)|\le \eta/{R_1}.
$$
Let us denote by $A_\vv$ the event under the probability sign in \ef{9.63}. Then
$$
I_2\le R_1 \pp(A_\vv)+\e(\ch_{A_\vv^c}\int_0^{R_1}|\psi_n(\uu(t;\vv))|\dd t).
$$
Combining \ef{9.63} with last two inequalities, we get $I_2\le 2\eta$, so that we have $I_1+I_2\le 3\eta$. Using this with \ef{9.64}, we see that
$$
\sup_{\vv\in B_R}|\elll\psi_n(\vv)|\le 3\eta\q\text{ for any }n\ge n_*(\eta).
$$
Since $\eta>0$ was arbitrary, the proof is complete.
\ep

\medskip
Let us note that in view of inequalities \ef{8.13} and \ef{9.22}, $\tilde\mu=\tilde\mu^\es$ satisfies $0<\tilde\mu(1)<\iin$.
We shall denote by $\hat\mu$ the normalization of $\tilde \mu$, that is
\be\label{9.68}
\hat\mu(\psi)=\tilde\mu(\psi)/\tilde\mu(1).
\ee 
Thanks to Lemma \ref{9.58}, $\hat\mu$ is continuous from $b_0(\h)$ to $\rr$.
For any Borel subset $\Gamma\subset \h$, we shall write
\be\label{9.66}
\hat\mu(\Gamma)=\hat\mu(\ch_\Gamma).
\ee
This notation will not lead to a confusion. 
\begin{lemma}\label{9.67}For any $\Gamma\subset \h$, we have
\be\label{9.65}
\mu(\dt\Gamma)\le\hat\mu(\dt\Gamma)\le\hat\mu(\bar\Gamma)\le\mu(\bar\Gamma),
\ee
where $\mu=\mu^\es$ is the stationary measure of the process $S^\es(\cdot)$.
\end{lemma}
\bp
We note that it is sufficient to show that for any closed set $F\subset\h$, we have
\be\label{9.71}
\hat\mu(F)\le\mu(F).
\ee

\medskip
{\it Step~1.}
Let us first show that for any bounded Lipschitz continuous function $\psi:\h\to\rr$, we have
\be\label{9.70}
\hat\mu(\psi)=(\psi, \mu).
\ee
Indeed, in view of \ef{9.61}, we have
\be\label{9.78'''}
\hat\mu(P_s\psi)=\hat\mu(\psi)\q\text{ for any }\psi\in b_0(\h).
\ee
In particular, this relation holds for any bounded Lipschitz function $\psi$ in $\h$. Moreover, it follows from inequality (1.3) in \cite{DM2014}, that for any such $\psi$, $P_s\psi$ converges to $(\psi,\mu)$ as $s\to\iin$ in the space $b_0(\h)$. Since $\hat\mu$ is continuous from $b_0(\h)$ to $\rr$, 
this implies 
$$
\hat\mu(\psi)=\hat\mu(P_s\psi)\to \hat\mu((\psi, \mu))=(\psi,\mu).
$$

\medskip
{\it Step~2.}
Now assume that inequality \ef{9.71} is not true, and let $F\subset \h$ closed and $\eta>0$ be such that
\be\label{9.72}
\hat\mu(F)\ge \mu(F)+\eta.
\ee
Let $\ch_F\le\psi_n\le 1$ be a sequence of Lipschitz continuous functions that converges pointwise to $\ch_F$ as $n\to \iin$. For example, one can take
$$
\psi_n(\uu)=\f{d_\h(\uu, F^c_{1/n})}{d_\h(\uu, F^c_{1/n})+d_\h(\uu,F)},
$$
where $F_r$ stands for the open $r$-neighborhood of $F$. Thanks to relation \ef{9.70}, inequality \ef{9.72} and monotonicity of $\hat\mu$, we have
$$
(\psi_n,\mu)=\hat\mu(\psi_n)\ge\hat\mu(\ch_F)=\hat\mu(F)\ge \mu(F)+\eta.
$$
However, this is impossible, since $(\psi_n, \mu)$ tends to $\mu(F)$ in view of the Lebesgue theorem on dominated convergence. The proof is complete.
\ep

\section{Proof of Proposition \ref{9.4}}
In view of Lemma \ref{9.67}, it is sufficient to prove that $\hat\mu=\hat\mu^\es$ satisfies
\be\label{9.78}
\exp(-(\vvv(\hat\uu_j)+\beta)/\es)\leq\hat\mu^\es(g_j)\leq \exp(-(\vvv(\hat\uu_j)-\beta)/\es).
\ee
\subsection{Upper bound} 
First note that by \eqref{9.59}-\ef{9.60}, for all $j\leq \ell+1$, we have
\be\label{8.2}
\tilde\mu^\es(\ch_{g_j})=\lm^\es(\e\int_0^{\tau_1^\es}\ch_{g_j}(S^{\es}(t;\cdot))\dd t)=\lm^\es(\ch_{\p g_j}(\cdot)\e\int_0^{\sigma_0^\es}\ch_{g_j}(S^{\es}(t;\cdot))\dd t).
\ee
In particular
\be\label{8.3}
\tilde\mu^\es(\ch_{g_j})\leq\lm^\es(\ch_{\p g_j}) \sup_{\vv\in\p g_j}\e_{\vv}\sigma_0^\es\leq\lm^\es(\ch_{\p g_j}) \sup_{\vv\in\tilde g}\e_{\vv}\sigma_0^\es.
\ee
On the other hand, we have
\be\label{8.4}
\tilde\mu^\es(1)\geq \tilde\mu^\es(\ch_{g'})\geq (1-\lm^\es(\ch_{\p g_{\ell+1}}))\min_{j\leq \ell} \inf_{\vv\in \bar g_j}\e\int_0^{\sigma_0^\es}\ch_{g_j}(S^{\es}(t;\vv))\dd t.
\ee
We recall that $g'$ stands for the set $g\backslash g_{\ell+1}$. We need the following result proved at the end of this section.
\begin{lemma}\label{8.21}
For any $\rho_*>0$ there exist $0<\rho'_1<\rho'_0<\rho_1<\rho_0<\rho_*$ such that for $\es<<1$ we have
\begin{align}
 \sup_{\vv\in\tilde g}\e_{\vv}\sigma_0^\es&\leq \exp(\beta/\es)\label{8.5},\\
\inf_{\vv\in \bar g_j}\e\int_0^{\sigma_0^\es}\ch_{g_j}(S^{\es}(t;\vv))\dd t&\geq \exp(-\beta/\es)\q\text{ for any } j\leq \ell.\label{8.6}
\end{align}
\end{lemma}
We first  note that $\vvv(\hat\uu_{\ell+1})$ is positive. Indeed, for any Lyapunov stable $\hat\uu_i$, the quantity $V(\hat\uu_i, \hat\uu_{\ell+1})$ is positive, and in view of \ef{8.45}, we have $\vvv(\hat\uu_{\ell+1})\ge \min V(\hat\uu_i, \hat\uu_{\ell+1})$, where the minimum is taken over $i\le \ell$ such that $\hat\uu_i$ is stable. Therefore, decreasing $\beta>0$ if needed, we may assume that $\vvv(\hat\uu_{\ell+1})\ge 2\beta$. In view of \ef{8.27'}, we have
\be\label{8.9}
\lm^\es(\ch_{\p g_{\ell+1}})\leq \exp (-(\vvv(\hat\uu_{\ell+1})-\beta)/\es)\leq\exp(-\beta/\es)\leq\f{1}{2}.
\ee
Combining this with inequalities \eqref{8.4} and \eqref{8.6}, we infer
\be\label{8.13}
\tilde\mu^\es(1)\geq \f{1}{2}\exp(-\beta/\es).
\ee
Further, using inequalities \ef{8.27'}, \eqref{8.3} and \eqref{8.5}, we get
\be\label{8.12}
\tilde\mu^\es(\ch_{g_j})\leq \exp(\beta/\es)\lm^\es(\ch_{\p g_j})\leq \exp(-(\vvv(\hat \uu_j)-2\beta)/\es).
\ee
Finally, combining this with \eqref{8.13}, we derive
$$
\hat\mu^\es(g_j)=\hat\mu^\es(\ch_{g_j})\leq 2\exp(\beta/\es) \exp(-(\vvv(\hat \uu_j)-2\beta)/\es)\leq \exp(-(\vvv(\hat \uu_j)-4\beta)/\es),
$$
where inequality holds for all $j\leq \ell+1$ and $\es<<1$.

\subsection{Lower bound}
We shall first establish the bound for $g_j$, $j\leq \ell$, and show that this implies the necessary bound for $g_{\ell+1}$. In view of \eqref{8.27'}, \eqref{8.2} and \eqref{8.6}, we have
\be\label{8.11}
\tilde\mu^\es(\ch_{g_j})\geq \exp(-(\vvv(\hat\uu_j)+2\beta)/\es).
\ee
On the other hand, by \eqref{8.12} we have
$$
\tilde\mu^\es(\ch_{g'})\leq \ell\exp(2\beta/\es).
$$
Note also that thanks to Lemmas \ref{9.47} and \ref{9.67}, we have that
$$
\hat\mu^\es(g')\ge\mu^\es(g')\geq \exp(-\beta/\es)\q\text{ for } \es<<1.
$$
Indeed, by definition, $g'$ represents the $\rho_1$-neighborhood of the set $\E=\{\hat\uu_1,\ldots, \hat\uu_\ell\}$.
It follows from the last two inequalities that
$$
\tilde\mu^\es(1)=\tilde\mu^\es(\ch_{g'})/\hat\mu^\es(g')\leq \ell\exp(3\beta/\es)\leq \exp(4\beta/\es)\q\text{ for } \es<<1.
$$
Finally, combing this inequality with \eqref{8.11}, we infer
\be\label{8.14}
\hat\mu^\es(g_j)\geq \exp(-(\vvv(\hat\uu_j)+6\beta)/\es),
\ee
where inequality holds for all $j\leq \ell$ and $\es<<1$.

\medskip
 We now show that inequality \eqref{8.14} implies 
\be\label{8.15}
\hat\mu^\es(g_{\ell+1})\geq \exp(-(\vvv(\hat\uu_{\ell+1})+8\beta)/\es)\q\text{ for }\es<<1.
\ee
We assume that $\vvv(\hat \uu_{\ell+1})<\iin$. First note that in view of \eqref{8.45}, we have 
\be\label{8.16}
\vvv(\hat \uu_{\ell+1})=\min_{i\leq \ell}[W_\ell(\hat \uu_i)+V(\hat\uu_i, \hat\uu_{\ell+1})]-\min_{i\leq \ell}W_\ell(\hat\uu_i).
\ee
Let us find $m\leq \ell$ such that
\be\label{8.17}
W_\ell(\hat \uu_m)+V(\hat\uu_m, \hat\uu_{\ell+1})=\min_{i\leq \ell}[W_\ell(\hat \uu_i)+V(\hat\uu_i, \hat\uu_{\ell+1})].
\ee

By definition of $V$, there is a finite time $T>0$ and an action $\ph\in L^2(0,T;H_\vartheta)$ such that
\be\label{5.3}
J_T(\ph)\leq V(\hat\uu_m, \hat\uu_{\ell+1})+\beta,\q |S^\ph(T;\hat\uu_m)-\hat\uu_{\ell+1}|_\h<\rho'_1/4.
\ee
Since the operator $S^\ph$ continuously depends on the initial point, there is $\kp>0$ such that 
$$
|S^\ph(T;\uu)-\hat\uu_{\ell+1}|_\h<\rho'_1/2,
$$
provided $|\uu-\hat\uu_m|_\h\leq\kp$. It follows from this inequality, relation \ef{9.78'''} and monotonicity of $\hat\mu^\es$ that
\begin{align}
\hat\mu^\es(g_{\ell+1})&=\hat\mu^\es(\ch_{g_{\ell+1}})=\hat\mu^\es(P_T\ch_{g_{\ell+1}})\ge \hat\mu^\es(\ch_{\bar B(\hat\uu_m, \kp)} P_T\ch_{g_{\ell+1}})\notag\\
&=\hat\mu^\es(\ch_{\bar B(\hat\uu_m, \kp)}(\cdot)\pp (|S^{\es}(T;\cdot)-\hat\uu_{\ell+1}|_\h<\rho'_1))\notag\\
&\geq  \hat\mu^\es(\ch_{\bar B(\hat\uu_m, \kp)}(\cdot)\pp(|S^{\es}(T;\cdot)-S^\ph(T;\cdot)|_\h<\rho'_1/2)\notag\\
&\ge\inf_{\uu\in \bar B(\hat\uu_m, \kp)}\pp(|S^{\es}(T;\uu)-S^\ph(T;\uu)|_\h<\rho'_1/2)\hat\mu^\es(B(\hat\uu_m, \kp))\label{5.6}.
\end{align}
In view of Theorem \ref{6.1} (applied to the set $B= \bar B(\hat\uu_m, \kp))$, we can find $\es_1=\es_1(\hat\uu_{\ell+1},\kp,\rho'_1,T)>0$ such that for all $\uu\in\bar B(\hat\uu_m, \kp)$, we have
$$
\pp(|S^{\es}(T;\uu)-S^\ph(T;\uu)|_\h<\rho'_1/2)\geq\exp(-(J_T(\ph)+\beta)/\es) \q\text {for }\es\leq \es_1.
$$
It follows that
\be\label{9.17}
\hat\mu^\es(g_{\ell+1})\geq\exp(-(J_T(\ph)+\beta)/\es) \hat\mu^\es(B(\hat\uu_m, \kp)) .
\ee
Combining this with first inequality of \eqref{5.3} and \eqref{5.6}, we get
$$
\hat\mu^\es(g_{\ell+1})\geq \hat\mu^\es  (B(\hat\uu_m, \kp))\exp(-(V(\hat\uu_m,\hat\uu_{\ell+1})+2\beta)/\es).
$$
Further, using this inequality and \eqref{8.14} \,\footnote{\, Recall that this inequality is true for any neighborhood of $\hat\uu_j,$ for $j\leq \ell$.} with $j=m$, we infer
$$
\hat\mu^\es(g_{\ell+1})\geq \exp(-(\vvv(\hat\uu_m)+V(\hat\uu_m,\hat\uu_{\ell+1})+8\beta)/\es).
$$
To complete the proof, it remains to note that thanks to \ef{8.17}, we have
\begin{align*}
\vvv(\hat\uu_m)+V(\hat\uu_m,\hat\uu_{\ell+1})&=W_\ell(\hat\uu_m)+V(\hat\uu_m,\hat\uu_{\ell+1})-\min_{i\leq \ell}W_\ell(\hat\uu_i)\\
&=\min_{i\leq \ell}[W_\ell(\hat \uu_i)+V(\hat\uu_i, \hat\uu_{\ell+1})]-\min_{i\leq \ell}W_\ell(\hat\uu_i)=\vvv(\hat\uu_{\ell+1}).\end{align*}

\subsection{Proof of Lemma \ref{8.21}} 
{\it Step~1: Derivation of \eqref{8.5}.} Let us fix $i\leq \ell+1$ and let $j\leq \ell$ be an integer different from $i$. Using Proposition \ref{1.88}, it is not difficult to show that $V(\hat\uu_i,\hat\uu_j)<\iin$ (see the derivation of \ef{9.12}). Therefore, if we denote by $d$ the distance between $\hat\uu_i$ and $\hat\uu_j$, we can find a finite time $T>0$ and an action $\ph$ defined on $[0, T]$ such that $S^\ph(T, \hat\uu_i)\in B(\hat\uu_j, d/2)$ and
$$
\int_0^T|\ph(s)|_{H_\vartheta}^2\dd s\leq V(\hat\uu_i,\hat\uu_j)+1.
$$
Let us find $t_*>0$ so small that 
\be\label{8.22}
\int_0^{t_*}|\ph(s)|_{H_\vartheta}^2\dd s\leq \beta\q\text{ and }\q \tilde \uu_i\neq\hat\uu_i,
\ee
where we set $\tilde \uu_i=S^\ph(t_*,\hat\uu_i)$. Further, let $\tilde\rho>0$ be so small, that $|\hat\uu_i-\tilde\uu_i|_\h\geq 4\tilde\rho$. And finally, let $0<\rho_*<\tilde\rho$ be such that
\be\label{8.23}
S^\ph(t_*;\uu)\in B(\tilde\uu_i,\tilde\rho)\q\text{ for any }\uu\in B(\hat\uu_i, \rho_*).
\ee
Now notice that if the trajectory $S^\es(t;\vv)$ issued from $\vv\in B(\hat\uu_i,\rho_*)$ is in the $\tilde\rho$-neighborhood of $S^\ph(t,\vv)$ in $C(0,t_*;\h)$ distance, then $S^\es(t_*,\vv)\notin B(\hat\uu_i,\rho_*)$. Moreover, if we denote by $\ph_*$ the restriction of $\ph$ on $[0, t_*]$, then for any $\vv$ in $B(\hat\uu_i,\rho_*)$, we have
$$
I_{t_*}(S^{\ph_*}(\cdot;\vv))=J_{t_*}({\ph_*})=\int_0^{t_*}|\ph(s)|_{H_\vartheta}^2\dd s\leq \beta.
$$
Applying Theorem \ref{6.1}, we derive
$$
\sup_{\vv\in B(\hat\uu_i,\rho_*)}\pp_\vv(\tau^\es_{\text{exit}}>t_*)\leq 1-\exp(-2\beta/\es),
$$
where we denote by $\tau^\es_{\text{exit}}(\vv)$ the time of the first exit of the process $S^\es(\cdot;\vv)$ from $B(\hat\uu_i,\rho_*)$. Now using the Markov property, we infer
\be\label{9.27}
\sup_{\vv\in B(\hat\uu_i,\rho_*)}\pp_\vv(\tau^\es_{\text{exit}}>n\,t_*)\leq (1-\exp(-2\beta/\es))^n,
\ee
which implies \eqref{8.5}.

 \bigskip
{\it Step~2: Proof of \eqref{8.6}.} Let us fix any stationary point $\hat\uu_i$ and let $\rho_*>0$. Given any $0<\rho_1<\rho_*$ let us find  $0<\rho_2<\rho_1$ such that for any $\vv\in \bar B(\hat\uu_i, \rho_2)$ we have $S(t;\vv)\in B(\hat\uu_i,\rho_1/2)$ for all $t\in [0,1]$. We assume that $\rho_*>0$ is so small that the conclusion of Lemma \ref{8.28} holds. We use the following construction: given any point $\vv\in \bar B(\hat\uu_i,\rho_1)$ we denote by $\tilde\ph_\vv$ the action defined on the time interval $[0, T+1]$ that coincides with $\ph_\vv$ on $[0,T]$ and vanishes on $[T, T+1]$. Then we have
\be\label{8.26}
I_{T+1}(S^{\tilde\ph_\vv}(\cdot;\vv))=J_{T+1}(\tilde\ph_\vv)=J_T(\ph_\vv)\leq\beta\q\text{ for any }\vv\in \bar B(\hat\uu_i,\rho_1).
\ee
Now let us take any $\rho_0\in (\rho_1, \rho_*)$ and let $\De<\min\{(\rho_0-\rho_1),\rho_2/2)\}$ be any positive number. Then by construction we have the following: if the trajectory $S^\es(t;\vv)$ is in the $\De$-neighborhood of $S^{\tilde\ph_\vv}(t;\vv)$ in the $C(0, T+1;\h)$ distance, then it remains in $\tilde g_i\equiv B(\hat\uu_i,\rho_0)$ for all $t\in [0,T+1]$ and moreover, it belongs to $g_i\equiv B(\hat\uu_i,\rho_1)$ for all $t\in[T,T+1]$. Therefore, we have
$$
\inf_{\vv\in \bar g_i}\e\int_0^{\sigma_0^\es}\ch_{g_j}(S^{\es}(t;\vv))\dd t\geq \pp(A), 
$$
where
$$
A=\{\om\in\Omega:\sup_{\vv\in \bar g_i} d_{C(0,T+1;\h)}(S^\es(\cdot;\vv), S^{\tilde\ph_\vv}(\cdot;\vv))<\De\}.
$$
Combining this with inequality \eqref{8.26} and Theorem \ref{6.1}, we arrive at \eqref{8.6}. The proof of Lemma \ref{8.21} is complete.

\section{A priori upper bound}\label{1.77}
In this section we establish Proposition \ref{9.18}. To simplify presentation, we first outline the main ideas.

\subsection{Scheme of the proof of Proposition \ref{9.18}}
\nt

\bigskip
{\it Compactness of level sets.} Let us suppose that we can prove the following: there is a constant $s\in (0, 1/2)$ such that
\be\label{1.75} 
|\uu_*|_{H^{s+1}(D)\times H^s(D)}\leq C(M)\q \text{ for any } \uu_*\in\{V_\aaa\leq M\},
\ee
where $H^s$ stands for the scale of Hilbert spaces associated with $-\de$.
Then the compactness of the embedding $H^{s+1}(D)\times H^s(D)\hookrightarrow\h$ implies that the level sets of $V_\aaa$ are relatively compact in $\h$. Thus, if inequality \eqref{1.75} is true, we only need to prove that the level sets  of $V_\aaa$ are closed.
Let $\uu_{*}^j$ be a sequence in $\{V_\aaa\leq M\}$ that converges to $\uu_*$ in $\h$. We need to show that $V_\aaa(\uu_*)\leq M$. By definition of $V_\aaa$, we have to show that for any positive constants $\eta$ and $\eta'$ there is an initial point $\uu_0\in\aaa$, a finite time $T=T_\eta>0$, and an action $\ph$ such that
\be\label{1.42}
J_{T}(\ph)\leq M+\eta'\q \text{ and }\q |S^{\ph}(T;\uu_{0})-\uu_*|_\h\leq \eta.
\ee
Let us fix $j$ so large that
\be\label{1.43}
|\uu_{*}^j-\uu_*|_\h\leq \eta/2.
\ee

Since $V_\aaa(\uu_*^j)\leq M$, there is a point $\uu_0\in\aaa$, a time $T=T_\eta>0$ and an action $\ph$ such that
$$
J_{T}(\ph)\leq M+\eta'\q \text{ and }\q |S^{\ph}(T;\uu_{0})-\uu_*^j|_\h\leq \eta/2.
$$
Combining this with inequality \eqref{1.43}, we derive \eqref{1.42}. The proof of inequality \eqref{1.75} relies on some estimates of the limiting equation and is carried out in the appendix.

\bigskip
{\it The bound \ef{9.19}}.
Due to the equivalence of \eqref{7.1} and \eqref{7.2}, we need to show that for any positive numbers $\De, \De'$ and $M$ there is $\es_*>0$ such that 
\be\label{1.99}
\mu^\es(\uu\in\h: d_{\h}(\uu,\{V_\aaa\leq M\})\geq\De)\leq\exp(-(M-\De')/\es)\q\text{ for } \es\leq\es_*.
\ee
From now on, we shall suppose that the constants $\De, \De'$ and $M$ are fixed.

\bigskip
{\it Reduction.} To prove \ef{1.99}, we first show that there is $\eta>0$ such that
\be\label{1.69}
\{\uu(t): \uu(0)\in \aaa_\eta, I_t(\uu_\cdot)\leq M-\De'\}\subset K_{\De/2}(M), \q t> 0,
\ee
where $\aaa_\eta$ stands for the open $\eta$-neighborhood of the set $\aaa$ and $K_\De(M)$ is  the open $\De$-neighborhood of the level set $\{V_\aaa\leq M\}$.
We then show that there is $R>0$ such that 
\be
\pP_1:=\mu^\es(B_R^c)\leq \exp(-M/\es)\q\text{ for }\es\le1.\label{1.30}
\ee
Once the constants $\eta$ and $R$ are fixed, we prove that there is $T_*>0$ such that
\be\label{2.16}
a:=\inf\{I_{T_*}(\uu_\cdot); \, \uu_\cdot\in C(0,T_*;\h), \, \uu(0)\in B_R, \,\uu(T_*)\in \aaa_\eta^c\}>0.
\ee

Taking inclusion \ef{1.69} and inequalities \ef{1.30}-\ef{2.16} for granted, let us show how to derive \ef{1.99}.

\bigskip
{\it Auxiliary construction.}
For any $n\geq 1$ introduce the set 
$$
E_n=\{\uu_\cdot\in C(0, nT_*;\h): \uu(0)\in B_R; \,\,\uu(kT_*)\in B_R\cap \aaa_\eta^c, \, k=1,\ldots, n\}.
$$
Let us mention that the idea of this construction is inspired by \cite{sowers-1992b} and $E_n$ is a modification of a set introduced by Sowers in that paper.
We claim that inequality \ef{1.69} and the structure of the set $E_n$ imply that for $n$ sufficiently large we have
\be\label{1.27}
\pP_2:=\sup_{\vv\in B_R}\pp(S^{\es}(\cdot;\vv)\in E_n)\leq\exp(-M/\es).
\ee

Indeed, in view of Theorem \ref{6.1}, to this end, it is sufficient to show that
\be\label{1.70}
\inf\{I_{nT_*}(\uu_\cdot); \, \uu_\cdot\in E_{n,T_*}\}> M.
\ee

We show that this inequality holds for any $n>(M+1)/a$. To see this, let us fix an integer $n$ satisfying this inequality and suppose that \eqref{1.70} is not true. Then there is an initial point $\vv\in B_R$ and an action $\ph$ defined on the interval $[0, nT_*]$ such that

$$
\f{1}{2}\int_{0}^{nT_*}|\ph(s)|_{H_\vartheta}^2\dd s<M+1
$$
and $S^\ph(jT_*;\vv)\in B_R\cap \aaa_\eta^c$ for all $j\in\{1,\ldots, n\}$. It follows from this inequality that there is $j\in\{1,\ldots, n-1\}$ such that
$$
\f{1}{2}\int_{jT_*}^{(j+1)T_*}|\ph(s)|_{H_\vartheta}^2\dd s<(M+1)/n<a.
$$
Therefore, the restriction of $\ph$ on the interval $[jT_*, (j+1)T_*]$ is an action whose energy is smaller than $a$ and that steers the point $\vv_1=S^\ph(jT_*;\vv)\in B_R$ to $\vv_2=S^\ph((j+1)T_*;\vv)\notin \aaa_\eta$. However, this is in contradiction with \eqref{2.16}. Inequality \ef{1.27} is thus established.

\bigskip
{\it Completion of the proof.} 
We now show that 
\be
\pP_3:=\int_{B_R}\pp(S^{\es}(t_*;\vv)\notin K_{\De}(M),\,S^{\es}(\cdot;\vv)\notin E_n)\mu^\es(d\vv)\leq \exp(-(M-2\De')/\es), \label{7.12}
\ee
where we set $t_*=(n+1)T_*$. Once this is proved, we will get \ef{1.99}. Indeed, by definition of the set $K_\De(M)$ and stationarity of $\mu^\es$, we have
\begin{align*}
\mu^\es(\uu\in\h: d_{\h}(\uu,\{V_\aaa\le M\})\geq\De)&=\mu^\es(\uu\in\h: \uu\notin K_\De(M))\\
&=\int_{\h}\pp(S^{\es}((n+1)T_*;\vv)\notin K_\De(M))\mu^\es(d\vv)\\
&\leq \pP_1+\pP_2+\pP_3.
\end{align*}
Combining inequalities \ef{1.30}, \eqref{1.27} and \eqref{7.12} we arrive at \ef{1.99}, where $\De'$
should be replaced by $3\De'$.

\medskip
To prove inequality \ef{7.12}, we first note that
\begin{align*}
\pP_3&\leq\int_{B_R}\pp\left(\bigcup_{j=1}^n \{S^{\es}(t_*;\vv)\notin K_{\De}(M),\,S^{\es}(jT_*;\vv)\in B_R^c\cup\aaa_\eta\}\right)\mu^\es(d\vv)\notag\\
&\leq \sum_{j=1}^n \int_{B_R}\pp(S^{\es}(t_*;\vv)\notin K_{\De}(M),\,S^{\es}(jT_*;\vv)\in B_R^c )\mu^\es(d\vv)\notag\\
&\q +\sum_{j=1}^n \int_{B_R}\pp(S^{\es}(t_*;\vv)\notin K_{\De}(M),\,S^{\es}(jT_*;\vv)\in\aaa_\eta)\mu^\es(d\vv):=\pP_3'+\pP_3''.
\end{align*}
By the stationarity of $\mu^\es$, the first term in the last inequality satisfies
\be\label{1.84}
\pP_3'\leq\sum_{j=1}^n\int_{\h}\pp(S^{\es}(jT_*;\vv)\in B_R^c )\mu^\es(d\vv)=n\,\pP_1\le n\exp(-M/\es).
\ee
Moreover, using inclusion \ef{1.69} and following \cite{sowers-1992b}, it is not difficult to prove (see Section \ref{9.40}) that
\be\label{9.8}
\pP_3''\leq\exp(-(M-\De')/\es)\q\text{ for }\es<<1
\ee
and thus to derive \ef{7.12}.

\bigskip
{\it Idea of the proof of \ef{1.69}-\ef{2.16}.} 
The proof of inequality \ef{1.30} is rather standard and relies on exponential estimates for solutions and a simple application of the Fatou lemma. 
The derivation of inclusion \eqref{1.69} is the most involved part in the proof. Without going into technicalities, we shall describe here the main ideas. We note that inclusion \eqref{1.69} clearly holds for $\eta=0$. Indeed, in this case $\uu(0)\in\aaa$, and since we have $I_t(\uu_\cdot)\leq M-\De'\leq M$, the point $\uu(t)$ is reached from the set $\aaa$ with action $\ph$ such that $J_t(\ph)\leq M$. It follows from the definition of $V_\aaa$ that $V_\aaa(\uu(t))\leq M$, so $d_\h(\uu(t),\{V_\aaa\le M\})=0$, and therefore we have \eqref{1.69}. So what we need to show is that if the initial point is sufficiently close to the attractor, then the inclusion \eqref{1.69} still holds. To prove this, we show that there is a flow $\hat \uu(t)$ issued from $\hat \uu(0)\in\aaa$ that remains in the $\De/2$-neighborhood of $\uu(t)$, and whose action function is $\De'$-close to that of $\uu(t)$. Once this is proved, the inclusion \eqref{1.69} will follow from the fact that $\hat \uu(t)\in \{V_\aaa\le M\}$, since it is reached from the set $\aaa$ at finite time $t$ with action whose energy is smaller than $M$. The construction of the flow $\hat \uu(t)$ relies on Proposition \ref{1.88}.

\medskip
As for the proof of inequality \eqref{2.16}, we first note that this inequality means the following: if we wait for sufficiently long time, then the energy needed to reach a  point outside $\eta$-neighborhood of the global attractor $\aaa$ is  positive uniformly with respect to the initial point in the ball $B_R$. The intuition behind this is that after sufficiently long time, the image of $B_R$ will be near the attractor $\aaa$, and the energy needed to steer the point close to the set $\aaa$ (say $\eta/2$-close) to a point outside its $\eta$-neighborhood, is positive. Let us finally mention that the fact that $V_\aaa$ vanishes only on the set $\aaa$ follows immediately from the definition of $V_\aaa$ and inequality \eqref{2.16}.

\subsection{Proof of inclusion \eqref{1.69}}\label{9.13}
\nt

{\it Step~1}. Let us suppose that \eqref{1.69} does not hold. Then there exist two sequences of positive numbers numbers $\eta_j\to 0$ and $T_j $, a sequence of initial points $(\uu^j_0)\subset\aaa_{\eta_j}$, and of action functions $(\ph^j)$ with $J_{T_j}(\ph^j)\leq M-\De'/2$, such that for each $j\geq 1$ the flow $\uu^j(t)=S^{\ph^j}(t;\uu_0^j)$ satisfies the inequality
\be\label{1.6}
d_\h(\uu^j(T_j),\{V_\aaa\le M\})\geq\De/2.
\ee

Let us also note that in view \eqref{1.38}, there is a positive constant $\mmm$ depending only on $\|h\|$ and $M$ such that for all $j\geq 1$ we have
\be\label{1.67}
\sup_{[0,T_j]} |\uu^j(t)|_\h\leq\mmm.
\ee

\bigskip
{\it Step~2}. For each $j\geq 1$, let us find $\vv^j_0\in\aaa$ such that $|\vv^j_0-\uu^j_0|_\h\leq\eta_j$ and introduce the intermediate flow $\vv^j(t)=[v(t),\dt v(t)]$ defined on the interval $[0,T_j]$ that solves 
\be\label{1.7}
\p_t^2 v+\gamma\p_t v-\de v+f(v)=h(x)+\ph^j+P_N[f(v)-f(u)], \q [v(0),\dt v(0)]=\vv^j_0,
\ee
where $N\geq 1$ is an integer to be chosen later and $u$ is the first component of $\uu^j(t)$. In view of Proposition \ref{1.88}, there is $N$ depending only on $\mmm$ such that for all $j\ge 1$ we have
\be\label{1.12}
|\vv^j(t)-\uu^j(t)|^2_\h\leq e^{-\al t}|\vv^j_0-\uu^j_0|_\h^2 \q\q\text{for all } t\in[0,T_j].
\ee

\bigskip
{\it Step~3}. Now let us fix $N=N(\mmm)$ such that we have \eqref{1.12}, and let us show that for $j>>1$  we have  
\be\label{1.16}
J_{T_j}(\hat\ph^j)\leq M, 
\ee
where we set 
\be\label{1.74}
\hat \ph^j=\ph^j+P_N[f(v)-f(u)].
\ee
By definition of $J$, we have
\be\label{1.14}
J_{T_j}(\hat\ph^j)=\f{1}{2}\int_0^{T_j}|\ph^j(s)+P_N[f(v(s))-f(u(s))]|_{H_\vartheta}^2\dd s.
\ee
We first note that
$$
|a+b|_{H_\vartheta}^2\leq p\,|a|_{H_\vartheta}^2+\f{p}{p-1}\,|b|_{H_\vartheta}^2,
$$
where $p>1$ is a constant to be chosen later. Therefore
\begin{align}
J_{T_j}(\hat\ph^j)&\leq \f{p}{2}\int_0^{T_j}|\ph^j(s)|_{H_\vartheta}^2\dd s+\f{p}{2(p-1)}\int_0^{T_j }|P_N[f(v(s))-f(u(s))]|_{H_\vartheta}^2\dd s\notag\\
&\leq p\,J_{T_j}(\ph^j)+C(N)\,\f{p}{p-1}\int_0^{T_j }|f(v(s))-f(u(s))|_{L^1}^2\dd s\label{1.15}.
\end{align}
By the H\"older and Sobolev inequalities, we have
$$
|f(v)-f(u)|_{L^1}^2\leq C_1\|u-v\|_1^2(\|u\|_1^2+\|v\|_1^2+1).
$$
Combining this with inequalities \eqref{1.67} and \eqref{1.12} we see that
$$
\int_0^{T_j }|f(v(s))-f(u(s))|_{L^1}^2\dd s\leq C_2 (\mmm^2+1)\int_0^{T_j}e^{-\al s}|\vv^j_0-\uu^j_0|_{\h}^2\dd s\leq C_3(\mmm^2+1)\eta_j^2.
$$
It follows from this inequality and \eqref{1.15}, and the fact that $N$ depends only on $\mmm$, that
$$
J_{T_j}(\hat\ph^j)\leq p\,J_{T_j}(\ph^j)+C(\mmm)\,\f{p}{p-1}\,\eta_j^2.
$$
Let us take 
$$
p=\f{M-\De'/4}{M-\De'/2}.
$$
Since $J_{T_j}(\ph^j)\leq M-\De'/2$, for $j$ large enough, we have $J_{T_j}(\hat\ph^j)\leq M$.

\bigskip
{\it Step~4}. We claim that for $j>>1$, we have
\be\label{1.17}
\vv^j(T_j)\in \{V_\aaa\le M\}.
\ee
Indeed, note that in view of \eqref{1.7}, we have $\vv^j(\cdot)=S^{\hat \ph^j}(\cdot; \vv_0^j)$. 
So the point $\vv^j(T_j)$ is reached from $\vv_0^j\in\aaa$ with action function $\hat \ph^j$ at finite time $T_j$. It follows from the definition of $V_\aaa$ and inequality \eqref{1.16} that $V_\aaa(\vv^j(T_j))\leq M$ for $j>>1$.

\bigskip
{\it Step~5}.
In view of inequality \eqref{1.12}, we have
$$
|\uu^j(T_j)-\vv ^j(T_j)|_\h\leq \eta_j \leq\De/4,
$$
provided $j\geq 1$ is large enough. Combining this with \eqref{1.17}, we see that
$$
d_\h(\uu^j(T_j), \{V_\aaa\le M\})\leq\De/4,
$$
which is in contradiction with \eqref{1.6}. The proof of  inclusion \eqref{1.69} is complete.

\subsection{Proof of inequality \eqref{2.16}}
Let us assume that \eqref{2.16} is not true, so for any $j\geq 1$ we have
$$
\inf\{I_{j}(\uu_\cdot); \, \uu_\cdot\in C(0,j;\h), \, \uu(0)\in B_R, \,\uu(j)\in \aaa_\eta^c\}=0.
$$
Then for each $j\geq 1$ there is an initial point $\uu^j_0\in B_R$ and an action $\ph^j$ defined on the interval $[0,j]$ with energy smaller than $e^{-j^2}$ such that the flow $\uu^j(t)=S^{\ph_j}(t;\uu^j_0)$ satisfies 
\be\label{2.4}
\uu^j(j)\notin\aaa_\eta.
\ee

For each $j\geq 1$, let $\vv^j(t)=S(t)\uu^j_0$. Using a priori bounds of the NLW equation it is not difficult to show (see Section \ref{9.29} for the proof) that
\be\label{9.41}
|\vv^j(t)-\uu^j(t)|_\h^2\leq C\int_0^t\|\ph^j(s)\|^2 \exp(C s)\dd s\q\text{ for } t\in [0, j].
\ee
Taking $t=j$ in this inequality and using $J_j(\ph^j)\le e^{-j^2}$, we get
\begin{align}
|\vv^j(j)-\uu^j(j)|_\h^2
&\leq C\exp(C j)\int_0^j \|\ph^j(s)\|^2ds\notag\\
&\leq C_3\exp(C j)J_j(\ph^j)\leq C_3\exp(-j^2+C j)\leq \eta^2/4\label{9.42},
\end{align}
provided $j$ is sufficiently large.
Combining this with \ef{2.4}, we see that for $j>>1$, we have
\be\label{9.30}
S(t)\uu^j_0=\vv^j(j)\notin\aaa_{\eta/2}.
\ee
Since $\aaa$ is the global attractor of the semigroup $S(t)$, we have
$$
\sup_{\uu_0\in B_R}d_\h(S(t)\uu_0, \aaa)\to 0\q\q\text{ as } t\to\iin.
$$
This is clearly in contradiction with \ef{9.30}. Inequality \ef{2.16} is established.
\subsection{Derivation of \ef{9.8}}\label{9.40}
We follow the argument presented in \cite{sowers-1992b}. Let us fix any $\vv\in B_R$ and $j\leq n$, and denote by $A$ the event $\{S^{\es}(t_*;\vv)\notin K_{\De}(M),\,S^{\es}(jT_*;\vv)\in\aaa_\eta\}$. Then, by the Markov property, we have
\begin{align*}
\pp(A)=\e[\e(\ch_A)|\fff^\es_{jT_*}]&=\e[\ch_{\bar \vv\in\aaa_\eta}\cdot \e(\ch_{S^{\es}(t_*-jT_*;\bar \vv)\notin K_{\De}(M)})]\notag\\
&\leq\sup_{\vv_0\in\aaa_\eta}\pp(S^{\es}(t_*-jT_*; \vv_0)\notin K_{\De}(M)),
\end{align*}
where $\fff^\es_t$ is the filtration corresponding to the process $S^{\es}(t;\vv)$ and we set $\bar \vv=S^{\es}(jT_*;\vv)$. It follows that
\be\label{3.2}
\pP_3''\leq \sum_{j=1}^n\sup_{\vv_0\in\aaa_\eta}\pp(S^{\es}(jT_*; \vv_0)\notin K_{\De}(M)).
\ee

\medskip
\nt
For any $l>0, M_1>0$ and $\vv\in\h$, introduce the level set
$$
K_{\vv,l}(M_1)=\{\uu(\cdot)\in C(0,l;\h); \uu(0)=\vv,\, I_l(\uu(\cdot))\leq M_1\}.
$$
Let us show that for any $l>0$ and $\vv_0\in\aaa_\eta$ we have
\be\label{3.3}
\{\om: S^{\es}(l;\vv_0)\notin K_{\De}(M)\}\subset \{\om: d_{C(0,l;\h)}(S^{\es}(\cdot;\vv_0), K_{\vv_0,l}(M-\De'))\geq \De/2\}.
\ee
Indeed, let us fix any $\om$ such that $S^{\es}(l;\vv_0,\om)\notin K_{\De}(M)$, and let $\uu_\cdot$ be any function that belongs to $K_{\vv_0,l}(M-\De')$. Then in view of inclusion \eqref{1.69} we have $\uu(l)\in K_{\De/2}(M)$, so that
$$
d_{C(0,l;\h)}(S^{\es}(\cdot;\vv_0,\om), \uu_\cdot)\geq |S^{\es}(l;\vv_0,\om)-z(l)|_\h\geq \De/2.
$$
Since $\uu_\cdot\in K_{\vv_0,l}(M-\De')$ was arbitrary, we conclude that inclusion \eqref{3.3} holds.
It follows from Theorem \ref{6.1} (applied to the time interval $[0, l]$ and the set $B=\aaa_\eta$) that there is $\es(l)=\es(l,\De,M,\eta)>0$ such that
\be\label{3.4}
\sup_{\vv_0\in\aaa_\eta}\pp(d_{C(0,l;\h)}(S^{\es}(\cdot;\vv_0), K_{\vv_0,l}(M-\De'))\geq \De/2)\leq \exp(-(M-2\De')/\es),\,\, \es\leq \es(l).
\ee
Let $\es_1(\De,M,T_*,n,\eta)=\min\{\es(T_*),\ldots \es(nT_*)\}$. Then in view of inequalities \eqref{3.2} and \eqref{3.4} we have
$$
\pP_3''\leq n\exp(-(M-2\De')/\es)\leq\exp(-(M-3\De')/\es) \q\text{ for }\es\leq\es_1.
$$
Inequality \ef{9.8} is established with $\De'$ replaced by $3\De'$.

\subsection{Proof of inequality \eqref{1.30}}
Let us show that for $R=R(M)$ sufficiently large and $\es_*=\es_*(M)>0$ small, we have \eqref{1.30}.  To this end, let us first show that the stationary solutions $\vv(t)$ of equation \eqref{0.1} satisfy 
\be\label{1.23}
\e\exp(\kp\,\ees (\vv(t))\leq Q(\es\,\BBB,\|h\|)\leq Q(\BBB,\|h\|),
\ee
for any $\kp\leq (\es\,\BBB)^{-1}\al/2$, where $Q$ and $\BBB$ are the quantities entering Proposition \ref{1.96}. Replacing $b_j$ by $b_j/\sqrt{\es}$, we see that it is sufficient to prove this inequality for $\es=1$.
Note that we cannot pass directly to the limit $t\to\iin$ in inequality \eqref{1.97}, since we first need to guarantee that $\e\exp(\kp\ees (\vv(0))$ is finite. This can be done by a simple application of the Fatou lemma. Indeed, for any $N\geq 1$, let $\psi_N(\uu)$ be the function that is equal to $\exp(\kp\ees (\uu))$ if $\ees(\uu)\leq N$, and to $\exp(\kp N)$ otherwise. Let us denote by $\mu$ the law of $\vv(t)$, and let $l$ be any positive number. Then using the stationarity of $\mu$ and inequality $\psi_N(\uu)\leq\exp(\kp\ees(\uu))$ we see that
\begin{align}
\int_\h \psi_N(\uu)\mu(d\uu)&=\int_\h\int_\h \psi_N(\uu')P_t(\uu, d\uu')\mu(d\uu)\notag\\
&\leq \int_{\ees(\uu)\leq l}\int_\h \exp(\kp\,\ees (\uu')) P_t(\uu, d\uu')\mu(d\uu)+\exp(\kp N)\mu(\ees(\uu)>l)\notag\\
&=\iI_1+\iI_2\label{1.80},
\end{align}
where $P_t$ stands for the transition function of the Markov process. Note that
$$
\iI_1\leq \sup_{\ees(\uu_0)\leq l}\e\exp(\kp\,\ees(\uu(t;\uu_0))\leq \exp(\kp \,l-\al\,t)+Q(\BBB,\|h\|),
$$
where $\uu(t;\uu_0)$ stands for the trajectory of \eqref{0.1} with $\es=1$ issued from $\uu_0$, and we used inequality \eqref{1.97}. Combining this with \eqref{1.80}, we obtain
$$
\int_\h \psi_N(\uu)\mu(d\uu)\leq \exp(\kp \, l-\al\,t)+\exp(\kp N)\mu(\ees(\uu)>l)+Q(\BBB,\|h\|).
$$
Passing to the limits $t\to\iin$ and then $l\to\iin$, and using the equivalence $\ees(\uu)\to\iin\Leftrightarrow |\uu|_\h\to\iin$, we get
$$
\int_\h \psi_N(\uu)\mu(d\uu)\leq Q(\BBB,\|h\|).
$$
Finally, letting $N$ go to infinity, and using Fatou's lemma, we derive \eqref{1.23}.

\medskip
We are now ready to establish \eqref{1.30}. Indeed, it follows from inequalities \eqref{1.56} that $\ees(\uu)\geq \f{1}{2}|\uu|_\h^2-2\,C> R^2/4$, provided $\uu\in B_R^c$ and $R^2\geq 8\,C$. Therefore, by the Chebyshev inequality, we have that
$$
\mu^\es(B_R^c)\leq \mu^\es(\ees(\uu)>R^2/4)\leq\exp(-\kp\, R^2/4)\int_{\h}\exp(\kp\,\ees(\uu))\,\mu^\es(d\uu).
$$
Now taking $\kp=(\es\,\BBB)^{-1}\,\al/2$ in this inequality, using \eqref{1.23} and supposing that $R$ is so large that $R\geq 16\,\BBB M\, \al^{-1}$, we obtain
$$
\mu^\es(B_R^c)\leq Q(\BBB,\|h\|) \exp(-2M/\es)\leq \exp(-M/\es),
$$
provided $\es>0$ is small. Inequality \eqref{1.30} is thus established.

\subsection{Lower bound with function $V_\aaa$ in the case of a unique equilibrium}\label{9.20}
The goal of this section is to show that in the case when equation \ef{1.51} possesses a unique equilibrium, the function $V_\aaa$ given by \ef{8.49}-\ef{9.6} provides also a lower bound for $(\mu^\es)$ and thus governs the LDP. The proof is almost direct and in this case there is no need to use the Freidlin-Wentzell theory.

\bigskip
So let $\hat\uu$ be the unique equilibrium of \ef{1.51}. It follows that the attractor $\aaa$ of the semigroup corresponding to \ef{1.51} is the singleton $\{\hat\uu\}$. Combining this with the fact that $(\mu^\es)$ is tight and any weak limit of this family is concentrated on $\aaa=\{\hat\uu\}$, we obtain 
\be\label{9.10}
\mu^\es\rightharpoonup\De_{\hat \uu}.
\ee 
We now use this convergence to establish the lower bound.
Due to the equivalence of \eqref{7.3} and \eqref{7.4}, we need to show that for any $\uu_*\in\h$ and any positive constants $\eta$ and $\eta'$, there is $\es_*>0$ such that we have \,\footnote{\,We write $V_{\hat\uu}$ instead of $V_\aaa$, since $\aaa=\{\hat\uu\}$.}
\be\label{5.2}
\mu^\es(B(\uu_*,\eta))\geq \exp(-(V_{\hat\uu}(\uu_*)+\eta')/\es) \q\text {for }\es\leq \es_*.
\ee
We assume $V_{\hat\uu}(\uu_*)<\iin$, since the opposite case is trivial. By definition of $V$, there is a finite time $T>0$ and an action $\ph\in L^2(0,T;H_\vartheta)$ such that
$$
J_T(\ph)\leq V_{\hat\uu}(\uu_*)+\eta'\q\text{ and }\q |S^\ph(T;\hat\uu)-\uu_*|_\h<\eta/4.
$$
Since the operator $S^\ph$ continuously depends on the initial point, there is $\kp>0$ such that 
$|S^\ph(T;\uu)-\uu_*|_\h<\eta/2$,
provided $|\uu-\hat\uu|_\h\leq\kp$. 
It follows that (see \ef{5.6}-\ef{9.17})
$$
\mu^\es(B(\uu_*,\eta))\geq \mu^\es(B(\hat\uu,\kp)) \exp(-(V_{\hat\uu}(\uu_*)+2\eta')/\es).
$$
Combining this inequality with convergence \ef{9.10} and using the portmanteau theorem, we infer
$$
\mu^\es(B(\uu_*,\eta))\geq C(\kp) \exp(-(V_{\hat\uu}(\uu_*)+2\eta')/\es)\geq  \exp(-(V_{\hat\uu}(\uu_*)+3\eta')/\es).
$$
Replacing $\eta'$ by $\eta'/3$, we arrive at \ef{5.2}.

\section{Appendix}

\subsection{Global attractor of the limiting equation}
In this section we recall some notions from the theory of attractors and an important result concerning the global attractor of the semigroup $S(t)$ generated by the flow of equation \eqref{1.51}.
\bi
\item {\it Equilibrium points}
\ei
We say that $\hat u\in \h$ is an equilibrium point for $S(t)$ if $S(t)\hat \uu=\hat \uu$ for all $t\geq 0$.
\bi
\item {\it Complete trajectory}
\ei
A curve $\uu(s)$ defined for $s\in\rr$ is called a complete trajectory of the semigroup $(S(t))_{t\geq 0}$ if 
\be\label{4.11}
S(t)\uu(s)=\uu(t+s)\q\text{ for all } s\in\rr \text{ and }t\in\rr_+.
\ee
\bi
\item {\it Heteroclinic orbits}
\ei
A heteroclinic orbit is a complete trajectory that joins two different equilibrium points, i.e., $\uu(t)$ is a heteroclinic orbit if it satisfies \eqref{4.11} and there exist two different equilibria $\hat \uu_1$ and $\hat \uu_2$, such that $\uu(-t)\to \hat \uu_1$ and $\uu(t)\to\hat \uu_2$ as $t\to\iin$.

\bi
\item {\it The global attractor}
\ei
\nt
The set $\aaa\subset \h$ is called the global attractor of the semigroup $(S(t))_{t\geq 0}$ if it has the following three properties:

1) $\aaa$ is compact in $\h$ ($\aaa\Subset \h$).

2) $\aaa$ is an attracting set for $(S(t))_{t\geq 0}$, that is
\be\label{4.9}
d_\h(S(t)B,\aaa)\to 0\q \text{ as } t\to\iin,
\ee
for any bounded set $B\subset\h$, where $d_\h(\cdot,\cdot)$ stands for the Hausdorff distance in $\h$.

3) $\aaa$ is strictly invariant under $(S(t))_{t\geq 0}$, that is
 \be\label{4.10}
S(t)\aaa=\aaa\q\text{ for all } t\geq 0.
\ee
The following result gives the description of the global attractor of the semigroup $S(t)$ corresponding to \eqref{1.51}. We assume that the nonlinear term $f$ satisfies \eqref{1.54}-\eqref{1.56}. We refer the reader to Theorem 2.1, Proposition 2.1 and Theorem 4.2 in Chapter 3 of \cite{BV1992} for the proof.
\begin{theorem}\label{Th-attractor}
The global attractor $\aaa$ of the semigroup $(S(t))_{t\geq 0}$ corresponding to \eqref{1.51} is a connected set that consists of equilibrium points of $(S(t))_{\geq 0}$ and joining them heteroclinic orbits. Moreover, the set $\aaa$ is bounded in the space $[H^2(D)\cap H^1_0(D)]\times H^1_0(D)$.
\end{theorem}

\subsection{Large deviations for solutions of the Cauchy problem}\label{6.0}
In this section we announce a version of large deviations principle for the family of Markov processes  generated by equation \eqref{0.1}. Its proof is rather standard, and relies on the contraction principle and the LDP for the Wiener processes.

\medskip

Let $B$ be any closed bounded subset of $\h$ and let $T$ be a positive number. We consider the Banach space $\mathsf{\yyy}_{B,T}$ of continuous functions $y(\cdot,\cdot):B\times[0,T]\to \h$ endowed with the norm of uniform convergence.

\begin{theorem}\label{6.1}
Let us assume that conditions \eqref{1.54}-\eqref{1.57} are fulfilled. Then $(S^{\es}(\cdot;\cdot), t\in [0, T], \vv\in B)_{\es>0}$ regarded as a family of random variables in $\yyy_{B,T}$ satisfies the LDP with rate function $I_T:\yyy_{B,T}\to [0,\iin]$ given by 
$$
I_T(y(\cdot,\cdot))=\f{1}{2}\int_0^T |\ph(s)|_{H_\vartheta}^2\dd s
$$
if there is $\ph\in L^2(0,T;H_\vartheta)$ such that $y(\cdot,\cdot)=S^{\ph}(\cdot;\cdot)$, and is equal to $\iin$ otherwise.
\end{theorem}
We refer the reader to the book \cite{DZ1992} and the paper \cite{CM-2010} for the proof of similar results.
Let us note that in the announced form, Theorem \ref{6.1} is slightly more general compared to the results from mentioned works, since they concern the case when the set $B$ is a singleton $B\equiv\{\uu_0\}$. However, recall that the LDP is derived by the application of a contraction principle to the continuous map
 $\gi: C(0,T; H^1_0(D))\to\yyy_{\uu_0,T}$ given by  
$$
\gi(q(\cdot))= y(\uu_0,\cdot),\q \text{ where }\, y(\uu_0,\cdot)=(S^{\dt q}(\uu_0;\cdot); t\in [0,T]).
$$
Using the boundedness of $B$, it is not difficult to show that the map $\tilde\gi$ from $ C(0,T; H^1_0(D))$ to $\yyy_{B,T}$ given by
$$
\tilde\gi(q(\cdot))= y(\cdot,\cdot),\q \text{ where }\, y(\cdot,\cdot)=(S^{\dt q}(\uu_0;\cdot); \uu_0\in B,  t\in [0,T])
$$
is also continuous. This allows to conclude.

\subsection{Lemma on large deviations}
\begin{lemma}\label{9.3}
Let $(\mathfrak{m}^\es)_{\es>0}$ be an exponentially tight family of probability measures on a Polish space $\zzz$ that possesses the following property: there is a good rate function $\mathfrak{I}$ on $\zzz$ such that for any $\beta>0$, $\rho_*>0$ and $z\in\zzz$ there are positive numbers $\tilde\rho<\rho_*$ and $\es_*$ such that
\begin{align*}
\mathfrak{m}^\es (B_\zzz(z,\tilde\rho))&\leq\exp(-(\mathfrak{I}(z)-\beta)/\es),\\
\mathfrak{m}^\es (\bar B_\zzz(z,\tilde\rho))&\ge\exp(-(\mathfrak{I}(z)+\beta)/\es) \q\text{ for }\es\leq\es_*.
\end{align*}
Then the family $(\mathfrak{m}^\es)_{\es>0}$ satisfies the LDP in $\zzz$ with rate function $\mathfrak{I}$.
\end{lemma}

\bp
We first note that in view of equivalence of \eqref{7.3} and \eqref{7.4}, we only need to establish the upper bound, that is inequality \eqref{7.1}. Moreover, since the family $(\mathfrak{m}^\es)_{\es>0}$ is exponentially tight, we can assume that $F\subset\zzz$ is compact (see Lemma 1.2.18 in \cite{DZ2000}). Now let us fix any $\beta>0$ and denote by $\tilde\rho(z)$ and $\es_*(z)$ the constants entering the hypotheses of the lemma. Clearly, we have
$$
F\subset\bigcup_{z\in F}B_\zzz(z,\tilde\rho(z)).
$$
Since $F\subset\zzz$ is compact, we can extract a finite cover 
$$
F\subset\bigcup_{i=1}^n B_\zzz(z_i,\tilde\rho(z_i)).
$$
It follows that
$$
\mathfrak{m}^\es(F)\leq n\exp(-(\inf_{z\in F}\mathfrak{I}(z)-\beta)/\es)\leq \exp(-(\inf_{z\in F}\mathfrak{I}(z)-2\beta)/\es),
$$
for $\es\leq\es_*(n, z_1,\ldots, z_n)$. We thus infer
$$
\limsup_{\es\to 0}\es\ln \mathfrak{m}^{\es}(F)\leq-\inf_{z\in F} \mathfrak{I}(z)+2\beta.
$$
Letting $\beta$ go to zero, we arrive at \eqref{7.1}.
\ep

\subsection{Proof of some assertions}\label{9.21}

\nt

\bigskip
{\it Genericity of finiteness of the set $\E$.} By genericity with respect to $h(x)$ we mean that $\E$ is finite for any $h(x)\in \cal{C}$, where $\cal{C}$ is a countable intersection of open dense sets (and therefore $\cal{C}$ is dense itself) in $H^1_0(D)$. This property is well known, and the proof relies on a simple application of the Morse-Smale theorem, see e.g., Chapter 9 in \cite{BV1992}. Here we would like to mention that there are also genericity results with respect to other parameters. Namely, it is known that in the case $h(x)\equiv 0$ and $f(0)=0$ the property of finiteness of the set $\E$ is generic with respect to the boundary $\p D$; we refer the reader to Theorem 3.1 in \cite{Saut1983}. Finally, let us mention that in the one-dimensional case, the genericity holds also with respect to the nonlinearity $f$, see \cite{Bru-Chow1984}.

\bigskip
{\it Exponential moments of Markov times $\sigma_k^\es$ and $\tau_k^\es$}. Let $R>0$ be so large that $\tilde g\subset B_R$.
 We claim that there is $\De(\es)>0$ such that we have
\be\label{9.22}
\sup_{\vv\in B_R}\e_\vv\exp(\De\sigma_0^\es)<\iin,\q \sup_{\vv\in B_R}\e_\vv\exp(\De\tau_1^\es)<\iin.
\ee
Indeed, in view of inequality \ef{9.27}, we have
$$
\sup_{\vv\in \tilde g}\e_\vv\exp(\De\sigma_0^\es)=\sup_{\vv\in \tilde g}\e_\vv\left(\sum_{n=0}^\iin \ch_{nt_*\le \sigma_0^\es< (n+1)t_*}e^{\De\sigma_0^\es}\right)\le e^{\De t_*}\sum_{n=0}^\iin (q e^{\De t_*})^n<\iin,
$$
where we set  $q=1-\exp(-\beta/\es)$ and choose $\De>0$ such that $q e^{\De t_*}< 1$. 
Thus, the first inequality in \ef{9.22} is established. By the strong Markov property, to prove the second one, it is sufficient to show that there is $\tilde\De\in (0, \De]$ such that 
$$
\sup_{\vv\in B_R}\e_\vv\exp(\tilde\De\tau^\es_g)<\iin, 
$$
where $\tau^\es_g(\vv)$ is the first instant when $S^\es(t;\vv)$ hits the set $\bar g$. The above relation follows \footnote{\,If the origin is among the equilibria, we use inequality (2.18) in the form announced in \cite{DM2014}. We note, however, that the latter is true for a neighborhood of any point, not only the origin, which allows to conclude in the general case.} from inequality (2.18) of  \cite{DM2014}, and we arrive at \ef{9.22} with $\De$ replaced by $\tilde \De$.

\subsection{Proof of inequality \ef{1.75}}

\nt

{\it Step~1.} Let $\uu_*\in\{V_\aaa\leq M\}$. By definition of $V_\aaa$, for any $j\geq 1$ there is an initial point $\uu_{0}^j\in\aaa$, a finite time $T_j>0$, and an action $\ph^j$ such that
\be\label{1.32}
J_{T_j}(\ph^j)\leq M+1\q \text{ and }\q |S^{\ph^j}(T_j;\uu_{0}^j)-\uu_*|_\h\leq 1/j.
\ee
In view of the second of these inequalities, in order to prove \eqref{1.75}, it is sufficient to show that
\be\label{1.44}
|S^{\ph^j}(T_j;\uu_{0}^j)|_{\h^s}\leq C(M) \q \text{ for all }j\geq 1,
\ee
where we set $\h^s=H^{s+1}(D)\times H^s(D)$.

\bigskip
{\it Step~2.} 
By definition of $S^\ph(t;\vv)$, we have $S^{\ph^j}(T_j;\uu_{0}^j)=\uu(T_j)$,
where $\uu(t)=[u(t),\dt u(t)]$ solves
\be\label{1.33}
\p_t^2 u+\gamma \p_t u-\de u+f(u)=h(x)+\ph^j(t,x),\q \uu(0)=\uu_{0}^j\q t\in[0,T_j].
\ee
We claim that
\be\label{1.45}
\ees(\uu(t))\leq C(\|h\|_1,M)\q\text{ for } t\in [0,T_j],
\ee
where $\ees(\uu)$ is given by \ef{1.95}.
Indeed, let us multiply equation \eqref{1.33} by $\dt u+\al u$ and integrate over $D$. Using some standard transformations and the dissipativity of $f$, we obtain 
$$
\p_t \ees(\uu(t))\leq -\al\ees(\uu(t))+C_1(\|h\|^2+\|\ph^j(t)\|^2) \q\q\q t\in[0,T_j].
$$
Applying the Gronwall lemma to this inequality, we get
\begin{align}
\ees(\uu(t))&\leq \ees(\uu(0))e^{-\al t}+C_1\int_0^t (\|h\|^2+\|\ph^j(s)\|^2)e^{-\al (t-s)}\dd s\notag\\
&\leq \ees(\uu(0))e^{-\al t}+C_1\al^{-1}\|h\|^2+C_1\int_0^{T_j}\|\ph^j(s)\|^2\dd s\notag\\
&\leq \ees(\uu(0))e^{-\al t}+C_1\al^{-1}\|h\|^2+C_2|\ph^j|_{L^2(0,T_j;H_\vartheta)}^2\leq C_3(\|h\|_1,M)\label{1.38},~~~~~~~
\end{align}
where we used first inequality of \eqref{1.32}, and the fact that since the initial point $\uu(0)=\uu_0^j$ belongs to the global attractor, its norm is bounded by constant depending on $\|h\|_1$ (see Theorem \ref{Th-attractor}).
Inequality \eqref{1.45} is thus established.

\bigskip
{\it Step~3.} 
We are now ready to prove \eqref{1.44}.
To this end, we split $u$ to the sum $u=v+z$, where $z$ solves
\be\label{1.37}
\p_t^2 z+\gamma \p_t z-\de z+f(u)=0,\q [z(0),\dt z(0)]=0\q t\in[0,T_j].
\ee
It is well known (e.g., see \cite{BV1992, Har85}) that inequality \eqref{1.45} implies  
\be\label{1.92}
|[z(T_j),\dt z(T_j)]|_{\h^s}\leq|[z(t),\dt z(t)]|_{C(0,T_j;\h^s)}\leq C_s(\|h\|_1,M)
\ee
for any $s<1-\rho/2$. 
Thus, it is sufficient to show that
\be\label{1.81}
|[v(t),\dt v(t)]|_{\h^{\f{1}{2}}}\leq C(\|h\|_1,M)\q\text{ for all } t\in[0,T_j].
\ee
Let us first note that in view of \eqref{1.33} and \eqref{1.37}, $v$ solves
\be\label{1.48}
\p_t^2 v+\gamma \p_t v-\de v=h(x)+\ph^j(t,x),\q [v(0),\dt v(0)]=\uu_0^j\q t\in[0,T_j].
\ee
Multiplying this equation with $-\de(\dt v+\al v)$, we obtain that for all $t\in [0,T_j]$
\begin{align}
\p_t |[v(t),\dt v(t)]|_{\h^{\f{1}{2}}}^2 &\leq -\al |[v(t),\dt v(t)]|_{\h^{\f{1}{2}}}^2+C_4(\|(-\de)^{\f{1}{2}} h\|^2+\|(-\de)^{\f{1}{2}}\ph^j(t)\|^2)\notag\\
&\leq -\al |[v(t),\dt v(t)]|_{\h^{\f{1}{2}}}^2+C_5(\|h\|_1^2+\|\ph^j(t)\|_{H_\vartheta}^2),\label{1.94}
\end{align}
where we used the fact that the space $H_{\vartheta}$ is continuously embedded in $H^1$, since 
$$
|\ph|_{\tilde H^1}^2=\sum_{j=1}^\iin\lm_j(\ph,e_j)^2=\sum_{j=1}^\iin\lm_j b_j^2(b_j^{-2}(\ph,e_j)^2)\leq \sup(\lm_j b_j^2)|\ph|_{H_{\vartheta}}^2
\leq\BBB_1 |\ph|_{H_{\vartheta}}^2.
$$
Applying the Gronwall lemma to inequality \eqref{1.94} and using first relation of \eqref{1.32} together with $\uu^j_0\in\aaa$, we derive \eqref{1.81}. Inequality \ef{1.75} is established.

\subsection{Some a priori estimates}\label{9.29}

\nt

\medskip
{\it Exponential moment of solutions.}
Let $\vv(t)$ be a solution of equation \eqref{0.1} with $\es=1$. We shall denote by $\ees:\h\to\rr$ the energy function given by
\be\label{1.95}
\ees(\uu)=|\uu|_\h^2+2\int_D F(u_1)\dd x \q\text { for } \uu=[u_1, u_2]\in\h.
\ee
The next result on the boundedness of exponential moment of $\ees(\vv(t))$ is taken from \cite{DM2014}.
\begin{proposition}\label{1.96}
Let conditions \ef{1.54}-\ef{1.57} be fulfilled. Then, we have
\be\label{1.97} 
\e\exp(\kp\,\ees(\vv(t))\leq\e\exp(\kp\,\ees (\vv(0))e^{-\al t}+Q(\BBB,\|h\|),
\ee
where inequality holds for any $\kp\leq (2\BBB)^{-1}\al$, and the function $Q(\cdot, \cdot)$ is increasing in both of its arguments. Here $\BBB$ stands for the sum $\sum b_j^2$ and $\al>0$ is the constant from \ef{9.11}.
\end{proposition}

\bigskip
{\it Feedback stabilization result.}
Let us consider functions $\uu(t)$ and $\vv(t)$ defined on the time interval $[0, T]$ that correspond, respectively, to the flows of equations \eqref{7.10} and
\be\label{1.98}
\p_t^2 v+\gamma \p_t v-\de v+f(v)+P_N [f(u)-f(v)]=h(x)+\ph(t,x).
\ee

We suppose that either $\ph(t,x)$ belongs to $L^2(0,T;L^2(D))$ or its primitive with respect to time belongs to $C(0,T;H^1_0(D))$. Let $\mmm$ be a positive constant such that 
\be\label{9.26}
|\uu(t)|_\h\vee |\vv(0)|_\h\leq \mmm \q\q \text{ for all } t\in [0,T].
\ee
The following result is a variation of Proposition 4.1 in \cite{DM2014}.
\begin{proposition}\label{1.88}
Under the conditions \ef{1.54}-\ef{1.56}, there is an integer $N_*$ depending only on $\mmm$ such that for all $N\geq N_*$ we have
\be\label{1.87}
|\vv(t)-\uu(t)|_\h^2\leq e^{-\al t}|\vv(0)-\uu(0)|_\h^2 \q\text{ for all } t\in [0, T].
\ee
\end{proposition}
\bp
We first show that there is an integer $N_1$ depending only on $\mmm$ such that for all $N\geq N_1$ we have
\be\label{9.25}
|\vv(t)|_\h\leq 4\mmm \q\q \text{ for all } t\in [0,T].
\ee
To this end, let us introduce
\be\label{1.9}
\tau=\inf\{t\in[0,T]: | \vv(t)|_{\h}>4\mmm\},
\ee
with convention that the infimum over the empty set is $\iin$. Inequality \eqref{9.25} will be proved if we show that there is $N_1=N_1(\mmm)$ such that $\tau=\iin$ for all $N\geq N_1$. Note that in view of \ef{9.26}, we have $\tau>0$. Moreover, by definition of $\tau$, we have
\be\label{1.11}
|\uu(t)|_\h \vee |\vv(t)|_\h\leq 4\mmm \q\q\q\text{ for all } t\in[0, \tau\wedge T].
\ee
It follows from Proposition 4.1 in \cite{DM2014} applied to the interval $[0,\tau\wedge T]$ that there is an integer $N_1$ depending only on $\mmm$ such that for all $N\geq N_1$ we have
\be 
|\vv(t)-\uu(t)|^2_\h\leq e^{-\al t}|\vv(0)-\uu(0)|_\h^2 \q\q\q\text{ for all } t\in[0,\tau\wedge T].
\ee
Therefore we have
\be
|\vv(\tau\wedge T)|_\h\leq |\uu(\tau\wedge T)|_\h+|\uu(0)|_\h+|\vv(0)|_\h\le 3\mmm.
\ee
Combining this with definition of $\tau$, we see that $\tau=\iin$, and thus inequality \eqref{9.25} is proved.
It follows that for $N\ge N_1$, we have
$$
|\uu(t)|_\h \vee |\vv(t)|_\h\leq 4\mmm \q\q\q\text{ for all } t\in[0, T].
$$ 
Once again using Proposition 4.1, but this time on the interval $[0,T]$, we see that there is $N_*\geq N_1$ such that for all $N\geq N_*$ we have \eqref{1.87}.
\ep

\bigskip
{\it Auxiliary estimates. Proof of \ef{9.41}.}
The standard argument shows (see the derivation of \eqref{1.38}) that there is a constant $\mmm$ depending only on $R$ and $\|h\|$ such that for all $j\geq 1$ we have
\be\label{2.6}
\sup_{[0, j]}|\uu^j(t)|_\h+\sup_{[0, j]}|\vv^j(t)|_\h\leq \mmm.
\ee
We shall write $\uu^j(t)=[u(t), \dt u(t)]$ and $\vv^j(t)=[v(t), \dt v(t)]$. Note that the difference $\uu^j(t)-\vv^j(t)$ corresponds to the flow of equation 
\be\label{2.3}
\p_t^2 z+\gamma\p_t z-\de z+f(v+z)-f(v)=\ph^j, \q [z(0),\dt z(0)]=0 \q\q  t\in [0,j].
\ee
Multiplying this equation by $\dt z+\al z$ and integrating over $D$, we obtain
\begin{align}
\p_t |[z(t),\dt z(t)]|_\h^2&\leq -\al |[z(t),\dt z(t)]|_\h^2+C(\|\ph^j(t)\|^2+\|f(v(t)+z(t))-f(v(t))\|^2)\notag\\
&\leq C(\|\ph^j(t)\|^2+\|f(v(t)+z(t))-f(v(t))\|^2).\label{2.7}
\end{align}
By the H\"older and Sobolev inequalities, we have
$$
\|f(v+z)-f(v)\|^2\leq C_1\|z\|_1^2(\|u\|_1^2+\|v\|_1^2+1)\leq C_2|[z,\dt z]|_\h^2(\|u\|_1^2+\|v\|_1^2+1).
$$
Combining this with inequalities \eqref{2.6} and \eqref{2.7}, we derive
$$
\p_t |[z(t),\dt z(t)]|_\h^2\leq C(\|\ph^j(t)\|^2+|[z(t),\dt z(t)]|_\h^2),
$$
where the constant $C$ depends only on $\mmm$.
Applying the comparison principle to this inequality, we see that for all $t\in[0,j]$, we have
$$
|[z(t),\dt z(t)]|_\h^2\leq C\int_0^t\|\ph^j(s)\|^2 \exp(C s)\dd s.
$$
Recalling the definition of $z$, we arrive at \ef{9.41}.

\subsection{Proof of Lemma \ref{1.201}}
We shall carry out the proof for the most involved case when $\uu_1$ and $\uu_2$ are the endpoints of the orbit $\ooo$. We need to show that for any positive constants $a$ and $\eta$, we have 
\be\label{9.28}
\mu^{\es_j}(B(\uu_2,\eta))\ge \exp(-a/\es)\q\q\text{ for } j>>1.
\ee
So, let $a$ and $\eta$ be fixed and let $\tilde\uu(\cdot)$ be a complete trajectory such that
$$
\tilde \uu(-t)\to\uu_1,\q \tilde \uu(t)\to\uu_2 \q\text{ as } t\to\iin.
$$
Let us find $T_1\geq 0$ and $T_2\geq 0$ so large that for $\tilde \uu_1=\tilde \uu(-T_1)$ and $\tilde \uu_2=\tilde \uu(T_2)$, we have
\be\label{5.8}
|\tilde \uu_1-\uu_1|_\h<\eta/8,\q |\tilde \uu_2-\uu_2|_\h<\eta/8.
\ee
Consider the flow $\uu'(t)=\tilde \uu(t-T_1)$. It corresponds to the flow of \eqref{1.51} issued from $\uu_1$.
Introduce the intermediate flow $\vv(t)$ corresponding to the solution of
$$
\p_t^2 v+\gamma \p_t v-\de v+f(v)=h(x)+\ph,\q [v(0),\dt v(0)]=\uu_1.
$$
where $\ph=P_N[f(v)-f(u')]$, and $u'$ is the first component of $\uu'$. In view of Proposition \ref{1.88}, there is an integer $N=N(\hat \uu_1,\|h\|)$ such that
\be\label{5.11}
|\vv(t)-\uu'(t)|_\h^2\leq e^{-\al t}\,|\uu_1-\tilde \uu_1|_\h^2\leq e^{-\al t}\,\eta^2/64,
\ee
where we used first inequality of \eqref{5.8}. In particular, for $T=T_1+T_2$, we have
\be\label{5.12}
|\vv(T)-\uu'(T)|_\h\leq\eta/8.
\ee
Now let us note that by construction we have
$$
\uu'(T)=\tilde \uu(T-T_1)=\tilde \uu(T_2)=\tilde \uu_2, \q \vv(T)=S^\ph(T;\uu_1).
$$
Combining this with second inequality of \eqref{5.8} and \eqref{5.12}, we obtain
$$
|S^\ph(T;\uu_1)- \uu_2|_\h<\eta/4.
$$
By continuity of $S^\ph$, there is $\kp>0$ such that
\be\label{5.13}
|S^\ph(T;\uu)-\uu_2|_\h<\eta/2 \q\text{ for } \,|\uu-\uu_1|_\h\leq\kp.
\ee
Moreover, it follows from inequality \eqref{5.11} that the action $\ph$ satisfies
\begin{align}
J_{T}(\ph)&=\f{1}{2}\int_0^{T}|\ph(s)|^2_{H_\vartheta}\dd s\leq C(N)\int_0^{T}|f(v)-f(u')|_{L^1}^2\dd s\notag\\
&\leq C_1\, C(N)\int_0^{T}(\|u'\|_1^2+\|v\|_1^2+1)\|v-u'\|_1^2\dd s\notag\\
&\leq C_2\,C(N) \int_0^{T}|\vv-\uu'|_\h^2 e^{-\al s}\dd s\leq C_3\, C(N)\, \eta^2\leq a\label{9.12},
\end{align}
provided $\eta$ is sufficiently small. Using inequality \eqref{5.13} and stationarity of $\mu^\es$ (see the derivation of \eqref{5.6}) we get
\begin{align*}
\mu^{\es_j}(B(\uu_2,\eta))&\geq \int_{|\uu-\uu_1|\leq\kp}\pp(|S^{\es_j}(T;\uu)-S^\ph(T;\uu)|<\eta/2)\,\mu^{\es_j}(d\uu)\\
&\ge \exp(-2a)\mu^{\es_j}(B(\uu_1,\kp))\q\q\text{ for } j>>1,
\end{align*}
where we used Theorem \ref{6.1} with inequality \eqref{9.12}. Moreover, since the point $\uu_1$ is stable with respect to $(\mu^{\es_j})$, we have 
$$
\mu^{\es_j}(B(\uu_1,\kp))\ge\exp(-a/\es_j)\q\q\text{ for } j>> 1.
$$
Combining last two inequalities, we arrive at \ef{9.28}, with $a$ replaced by $3a$. Lemma \ref{1.201} is established.

\addcontentsline{toc}{section}{Bibliography}
\bibliographystyle{plain}
\def\cprime{$'$} \def\cprime{$'$}
  \def\polhk#1{\setbox0=\hbox{#1}{\ooalign{\hidewidth
  \lower1.5ex\hbox{`}\hidewidth\crcr\unhbox0}}}
  \def\polhk#1{\setbox0=\hbox{#1}{\ooalign{\hidewidth
  \lower1.5ex\hbox{`}\hidewidth\crcr\unhbox0}}}
  \def\polhk#1{\setbox0=\hbox{#1}{\ooalign{\hidewidth
  \lower1.5ex\hbox{`}\hidewidth\crcr\unhbox0}}} \def\cprime{$'$}
  \def\polhk#1{\setbox0=\hbox{#1}{\ooalign{\hidewidth
  \lower1.5ex\hbox{`}\hidewidth\crcr\unhbox0}}} \def\cprime{$'$}
  \def\cprime{$'$} \def\cprime{$'$} \def\cprime{$'$}

\end{document}